\documentclass[review,hidelinks,onefignum,onetabnum]{siamart220329}

\usepackage{lipsum}
\usepackage{amsfonts}
\usepackage{graphicx}
\usepackage{epstopdf}
\usepackage{algorithmic}
\ifpdf
  \DeclareGraphicsExtensions{.eps,.pdf,.png,.jpg}
\else
  \DeclareGraphicsExtensions{.eps}
\fi

\newsiamremark{remark}{Remark}
\newsiamremark{hypothesis}{Hypothesis}
\crefname{hypothesis}{Hypothesis}{Hypotheses}
\newsiamthm{claim}{Claim}

\headers{Stability theory of TASE-RK methods with inexact Jacobian}{D. Conte, J. Martin-Vaquero, G. Pagano and B. Paternoster}

\title{Stability theory of TASE-Runge-Kutta methods with inexact Jacobian\thanks{Submitted to the editors DATE.
\funding{GNCS-INdAM project; PRIN PNRR 2022 project (No. P20228C2PP); MUR PRIN 2020 project (No. 2020JLWP23); %“Integrated Mathematical Approaches to Socio–Epidemiological Dynamics”; %(CUP: E15F21005420006) 
MCIN European Union NextGenerationEU (PRTRC17.I1) project; MCIN Consejer\'ia de Educación, Junta de Castilla y Le\'on, QCAYLE project; Fundaci\'on Sol\'orzano, FS/5-2022 project; “Visiting Professors mobility Program” of the University of Salerno.}}}

\author{Dajana Conte\thanks{Department of Mathematics, University of Salerno, 84084, Fisciano, Italy
  (\email{dajconte@unisa.it}, \email{gpagano@unisa.it}, \email{beapat@unisa.it}).}
\and Jesus Martin-Vaquero\thanks{Department of Applied Mathematics and IUFFyM, University of Salamanca, E37700, Bejar, Spain
  (\email{jesmarva@usal.es}).}
\and Giovanni Pagano\footnotemark[2]
\and Beatrice Paternoster\footnotemark[2]}

\usepackage{amsopn}

\ifpdf
\hypersetup{
	pdftitle={Stability theory of TASE-Runge-Kutta methods with inexact Jacobian},
	pdfauthor={D. Conte, J. Martin-Vaquero, G. Pagano and B. Paternoster}
}
\fi

\usepackage{comment}
\usepackage{upgreek}

\newtheorem{example}{Test}

\makeatletter
\renewcommand*\env@matrix[1][\arraystretch]{%
	\edef\arraystretch{#1}%
	\hskip -\arraycolsep
	\let\@ifnextchar\new@ifnextchar
	\array{*\c@MaxMatrixCols c}}
\makeatother

\def\<#1>{\mathinner{\langle#1\rangle}}

\begin{document}

\nolinenumbers
	
	\maketitle
	
	% REQUIRED
	\begin{abstract}
		This paper analyzes the stability of the class of Time-Accurate and Highly-Stable Explicit Runge-Kutta (TASE-RK) methods, introduced in 2021 by Bassenne et al. (J. Comput. Phys.) for the numerical solution of stiff Initial Value Problems (IVPs). Such numerical methods are easy to implement and require the solution of a limited number of linear systems per step, whose coefficient matrices involve the exact Jacobian $J$ of the problem. To significantly reduce the computational cost of TASE-RK methods without altering their consistency properties, it is possible to replace $J$ with a matrix $A$ (not necessarily tied to $J$) in their formulation, for instance fixed for a certain number of consecutive steps or even constant. However, the stability properties of TASE-RK methods strongly depend on this choice, and so far have been studied assuming $A=J$.
		
		In this manuscript, we theoretically investigate the conditional and unconditional stability of TASE-RK methods by considering arbitrary $A$. To this end, we first split the Jacobian as $J=A+B$. Then, through the use of stability diagrams and their connections with the field of values, we analyze both the case in which $A$ and $B$ are simultaneously diagonalizable and not. Numerical experiments, conducted on Partial Differential Equations (PDEs) arising from applications, show the correctness and utility of the theoretical results derived in the paper, as well as the good stability and efficiency of TASE-RK methods when $A$ is suitably chosen.
	\end{abstract}
	
	% REQUIRED
	\begin{keywords}
		Runge-Kutta methods,
		TASE operators,
		W-methods,
		stiff problems,
		stability with inexact Jacobian,
		conditional stability,
		unconditional stability,
		stability diagrams,
		field of values. %\LaTeX
	\end{keywords}
	
	% REQUIRED
	\begin{MSCcodes}
		65L04, 65L06, 65L07, 65L20, 65M12, 65M20.
	\end{MSCcodes}
	
	\section{Introduction}\label{sec1}
	Models based on Ordinary Differential Equations (ODEs) and PDEs arising from applications are often characterized by severe stiffness. When solving stiff equations, the use of explicit methods in time is discouraged, as they need decidedly small step-sizes to avoid instability issues. The methods characterized by the best stability properties are the implicit ones, which however have a high computational cost as they require the solution of a system, of the size of the problem, of nonlinear equations at each time step. For implicit methods, several techniques have been introduced to practically work with systems of nonlinear equations of smaller dimensions, and to minimize the number of steps of the iterative process which is chosen to solve them, see e.g. \cite[Ch. IV.8, pp. 118--122]{HairerWanner}. However, such methods are usually not efficient for the solution of huge problems arising from applications.
	
	For this reason, for the solution of stiff problems it is customary to use, e.g., Implicit-Explicit (ImEx) \cite[Ch. IV.4]{Hundsdorferbook} or linearly implicit methods. Roughly speaking, ImEx methods require a splitting of the problem into a stiff part and a non-stiff part, addressing the former through an implicit scheme and the latter through an explicit scheme. Linearly implicit methods, on the other hand, require the solution of a fixed (small) number of linear systems per step, whose coefficient matrices involve the exact Jacobian of the problem or an approximation thereof.
	Well known and efficient linearly implicit schemes are the Rosenbrock \cite[Sec. II.1.5]{Hundsdorferbook} and W-methods \cite{Steihaug}.
	
	In this paper, we consider the class of TASE-RK methods, which has been recently introduced in \cite{Bassenne2021} and then improved in \cite{Calvo2021}. TASE-RK methods constitute a subclass of the W-methods \cite{Gonzalez-Pinto2023129}, and are particularly practical and simple to implement, as they have been designed for engineering applications. Interest in these numerical methods is growing, as evidenced by the large amount of related scientific manuscripts that have been published in the last years \cite{Aceto2023,Bassenne2021,Calvo2023,Calvo2021,Conte20231,Gonzalez-Pinto2023129}. Indeed, TASE-RK methods are quite cheap, flexible, and characterized by good accuracy and stability properties. In the manuscripts \cite{Aceto2023,Gonzalez-Pinto2023129}, it has been shown that TASE-RK schemes and their variants are competitive with some implicit RK, Rosenbrock and W-methods.
	
	The linear systems to be solved at each step in TASE-RK methods are of the form $(I-\omega k A) \vec{U} = \vec{F}$, where $\omega$ and $k$ are positive parameters, $\vec{U}$ is a vector to be determined and $\vec{F}$ is a known function; by definition, the matrix $A$ corresponds to the exact Jacobian $J$ of the problem. However, since the consistency of these numerical schemes is independent of the choice of $A$, TASE-RK methods are generally applied by fixing an arbitrary matrix in place of the Jacobian, preferably constant to drastically reduce their cost. In this way, the order of convergence is naturally preserved. Nevertheless, if this arbitrary matrix is not a good approximation of the exact Jacobian, there is a concrete risk of running into instability issues. Indeed, until now, the good stability properties of TASE-RK methods have been proved by assuming $A=J$. This motivates us to theoretically study the stability of TASE-RK methods when $A$ does not correspond to the Jacobian.

	\subsection{Problem, recent contributions, and novelty}
	The differential equations we are interested in are IVPs of the type
	\begin{equation}\label{ODEs}
	\begin{cases}
	\displaystyle \vec{u}_t= \vec{f}(t,\vec{u}), \\
	\displaystyle \vec{u}(t_0)=\vec{u}_0,
	\end{cases}
	\qquad \vec{f}:\mathbb{R} \times \mathbb{R}^d \rightarrow \mathbb{R}^d, \quad t \in [t_0,t_e].
	\end{equation}
	Problems of the form \eqref{ODEs} can also arise, of course, from the spatial semi-discretization of PDEs performed by means of several techniques, e.g. via finite differences or finite elements. There are many real-world phenomena that give rise to semi-discretized problems equipped with severe stiffness, such as chemical reactions, corrosion of materials \cite{Frasca-Caccia2023}, vegetation growth \cite{Conte2023}, and so on.
	
	To study the stability of TASE-RK methods we consider the test problem
	\begin{equation}\label{test_ODEs}
	\vec{u}_t= J\vec{u} + \vec{g}(t), \quad \vec{u}(t_0)=\vec{u}_0,
	\qquad \text{with} \quad J \in \mathbb{R}^{d \times d}, \quad  \vec{g}:\mathbb{R} \rightarrow \mathbb{R}^d, \quad t \in [t_0,t_e].
	\end{equation}
	As usual, we assume that all the eigenvalues of the matrix $J$ have negative real part and that the boundedness of the solution does not depend on the forcing term $\vec{g}$. This approach is standard for the study of the stability properties of numerical methods \cite{Butcherbook}; the matrix $J$ in \eqref{test_ODEs} plays the role of the Jacobian of the function $\vec{f}$ in \eqref{ODEs}.
	
	When $\vec{f}$ is characterized by a natural splitting into stiff and non-stiff parts, it is convenient to use numerical schemes such as ImEx or W-methods, which require the solution of linear systems at each step whose coefficient matrices do not necessarily involve the exact Jacobian of the problem. To analyze the stability of ImEx methods, it is customary to consider the test problem \eqref{test_ODEs} with $J=A+B$; the linear systems to be solved involve only the matrix $A$, which is linked to the stiff part of the vector field $\vec{f}$ in \eqref{ODEs}. For simplicity, $A$ and $B$ are generally assumed to be simultaneously diagonalizable, see e.g. \cite[Sec. 2.6]{Boscarino2015}. In recent years, however, Rosales et al. \cite{Rosales20172336} and Seibold et al. \cite{Seibold2019} have considered a class of ImEx Linear Multistep Methods (LMMs), studying their unconditional (linear) stability by removing the hypothesis of simultaneous diagonalizability of $A$ and $B$, exploiting connections between the field of values of a matrix \cite[Ch. I]{HornJohnsonFOV} and the stability diagrams of the schemes.
	
	In this manuscript, we derive the stability diagrams of the TASE-RK methods and prove that, even for them, it is possible to establish unconditional stability conditions which are linked to the field of values. Furthermore, we carry out a study of the conditional stability of TASE-RK schemes, both by avoiding the hypothesis of simultaneous diagonalizability of $A$ and $B$, and by considering it (since this leads to less restrictive conditions). Finally, we also comment on nonlinear stability. We underline that TASE-RK methods constitute a subclass of W-methods, for which stability is usually studied considering the Rosenbrock scheme associated with them (i.e. assuming that the exact Jacobian is used in the formulation), see e.g. \cite[Sec. 4]{LangNumerMath2013}. Therefore, our results could prove very useful and constitute a starting point for formally analyzing the stability of general W-methods. Indeed, to our knowledge, this is the first manuscript where the stability of a subclass of W-methods is in-depth studied with an inexact Jacobian, using the techniques that we employ.

	\subsection{Outline of the paper}
	
	This paper is organized as follows. In Section \ref{sec2} we recall the formulation and properties of TASE-RK methods. In Section \ref{sec3} we define the so-called conditional and unconditional stability diagrams related to TASE-RK methods, proving for them some geometric properties that will be very useful for the results derived subsequently. In Section \ref{sec4} we prove some results about the stability of TASE-RK methods in the hypothesis in which the matrices $A$ and $B$ are simultaneously diagonalizable. In Section \ref{sec5} we consider the general case, by removing the assumption of simultaneous diagonalizability of $A$ and $B$. In Section \ref{sec6} we give some hints on nonlinear stability and provide an example of application of the theory. In Section \ref{sec7} we perform numerical experiments, using the theoretical results derived in the paper to solve PDEs from applications, by means of TASE-RK methods with inexact Jacobian; the efficiency of these schemes will be confirmed by comparison with other linearly implicit numerical methods. In Section \ref{sec8} we briefly summarize the obtained results, discussing possible future research scenarios.

	\section{TASE-Runge-Kutta methods}\label{sec2}
	From now on, we consider the homogeneous discrete grid $\{t_n=t_0+nk; \ n=0,\ldots,N; \ t_N=t_e\}$ for the time integration of the IVP \eqref{ODEs}. In \cite{Bassenne2021}, the authors have designed the so-called TASE-RK methods for the solution of stiff IVPs, deriving from the combination of the following two steps:
	\begin{enumerate}
		\item the problem \eqref{ODEs} is modified premultiplying the vector field $\vec{f}$ by a suitably chosen matrix operator $T_p$, called TASE operator, which approximates the identity mapping up to a given order $p$;
		\item the modified problem is solved through an explicit RK method of order $p$.
	\end{enumerate}
	TASE operators are functions of the product $kJ$, where $k$ is the time step and $J$ is the Jacobian of $\vec{f}$ in \eqref{ODEs}, and of a vector $\vec{\omega}$ of free parameters to be set appropriately. Therefore, more rigorously, the first step means that a TASE operator fulfils
	\begin{equation}\label{Tp_identity}
	T_p(kJ;{\vec{\omega}})=I+O(k^p),
	\end{equation}
	where, from now on, $I$ stands for the identity matrix of the appropriate size. The operator $T_p$ must also be chosen in such a way that the TASE-RK method which is obtained by combining the two steps described above has good stability properties.
	
	In \cite{Bassenne2021}, the first family of TASE operators has been derived, through Richardson extrapolation starting from $T_1(kJ;\omega)=(I-\omega k J)^{-1}=I+O(k), \ \omega \in \mathbb{R}$ (free parameter), in order to satisfy \eqref{Tp_identity}. In this paper, we focus on the following family of TASE operators, which has been proposed in \cite{Calvo2021} to generalize the one derived in \cite{Bassenne2021}:
	\begin{equation}\label{TASE_Montijano}
	\begin{split}
	& T_p(kJ;{\vec{\omega}})=\sum_{j=1}^{p}\beta_j(I-\omega_j k J)^{-1}, \quad {\vec{\omega}}=(\omega_j) \in \mathbb{R} ^ p, \quad \omega_j > 0 \quad \forall j, \\ & \beta_j:=\Big(\frac{1}{\omega_j}\Big)^{p-1}/\prod_{l\neq j}\Big(\frac{1}{\omega_j}-\frac{1}{\omega_l}\Big), \quad j=1,\ldots,p.
	\end{split}
	\end{equation}
	The values of the coefficients $ \beta_j $ are a-priori fixed in such a way that $T_p$ satisfies \eqref{Tp_identity}, while the $\omega_j$ parameters are free and should be set for stability reasons.
	Other families of TASE operators have been proposed in
	\cite{Aceto2023,Calvo2023,Conte20231}.
	
	TASE-RK methods are widely used in engineering applications as they do not require implementation efforts: indeed, it is sufficient to know how to construct an explicit RK method (if it then gives stability issues, the problem is modified by introducing $T_p$).
	For ease of analysis and to avoid modifying the problem to be solved, we formulate $s$-stage TASE-RK methods for \eqref{ODEs} as
	\begin{equation}\label{TASE_RK}
	\begin{cases}
	\begin{split}
	\vec{U}_{n,i}=&\vec{u}_n+k\sum_{j=1}^{i-1} \alpha_{ij}T_{p}(kA;{\vec{\omega}})\vec{f}(t_n+c_jk,\vec{U}_{n,j}), \qquad i=1,\ldots,s,  \\
	\vec{u}_{n+1}=&\vec{u}_n+k\sum_{j=1}^{s}b_jT_{p}(kA;{\vec{\omega}})\vec{f}(t_n+c_jk,\vec{U}_{n,j}),
	\end{split}
	\end{cases}
	\end{equation}
	where $\vec{u}_n \approx \vec{u}(t_n)$ and $\vec{U}_{n,i} \approx \vec{u}(t_n + c_ik)$.
	This reformulation is much more natural than considering TASE-RK methods as the combination of the two steps described at the beginning of this section. Furthermore, it allows to simplify the analysis of the consistency and stability properties of the TASE-RK schemes, and to further highlight their link with the W-methods.
	
	\begin{remark}
		Note that \eqref{Tp_identity} ensures that the TASE-RK method preserves the order $p$ of the underlying explicit RK method. Exploiting this,
		we here consider a generalization of TASE-RK methods \eqref{TASE_RK},  replacing the Jacobian $J$ in $T_p$ with an arbitrary matrix $A$.
		Setting indeed $\beta_j$ as in \eqref{TASE_Montijano}, property \eqref{Tp_identity} continues to hold even with this replacement. Clearly, this generalization gives greater flexibility to the TASE-RK methods, and could allow a significant reduction in their computational cost if $A$ is chosen constant for each step, or at least for several consecutive steps. However, the study of stability becomes complicated when $A \neq J$.
	\end{remark}
	
	Up to now, the stability of TASE-RK methods \eqref{TASE_RK} has been investigated under the hypothesis that $A$ corresponds to the Jacobian of the problem. That is, considering the test equation $u_t=\eta u$, with $\eta \in \mathbb{C}$, Re$(\eta)\leq 0$, it has been assumed that $A=\eta$. With this hypothesis, the stability function of $s$-stage TASE-RK schemes of order $p=s\leq 4$ is given by \cite[p. 3]{Calvo2021}
	\begin{equation}\label{Stab_TASE}
	RT_p(z;{\vec{\omega}})=1+T_p(z;{\vec{\omega}})z+\frac{1}{2}(T_p(z;{\vec{\omega}})z)^2+\ldots+\frac{1}{p!}(T_p(z;{\vec{\omega}})z)^p, \quad \text{where} \ z=k\eta.
	\end{equation}
	Note that $RT_p$ depends on the free coefficients $\omega_j$ of the TASE operator $T_p$.
	Values of ${\vec{\omega}}$ which guarantee good stability when $RT_p$ is as in \eqref{Stab_TASE} have been determined in \cite{Calvo2021} (see Table \ref{tableTRK}). Recall that, for TASE-RK methods: \textit{A-stability} means $|RT_p(z)|\leq 1$ $\forall \ z\in \mathbb{C}^-$; \textit{strong A-stability} means A-stability with $|RT_p(\infty)|:=\lim_{z\rightarrow - \infty} |RT_p(z)|<1$ ($z \in \mathbb{R}$); \textit{L-stability} means strong A-stability with $|RT_p(\infty)|=0$; \textit{A$(\theta)$-stability} means $|RT_p(z)|\leq 1$ $\forall \ z\in S_\theta$, $S_\theta:=\{z\in \mathbb{C}  : |\text{arg}(-z)|\leq \theta \}$; \textit{strong A$(\theta)$-stability} and \textit{L$(\theta)$-stability} follow from strong A-stability, L-stability, and A$(\theta)$-stability.
	
	\begin{table}[h!]
		\centering
		\caption{Stability properties of $s$-stage TASE-RK methods of order $ p=s \leq 4 $, setting ${\vec{\omega}}$ as in \cite{Calvo2021}.} \label{tableTRK}
		\begin{tabular}{c ||  c | c | c || c}
			\hline
			$p$ &  Stability properties & $\theta$ & $RT_p(\infty)$ & ${\vec{\omega}}$ \\
			\hline
			2   & Strong A-stability & $90^o$ & 0.5 &  $(3,1.5)$ \\
			3  & L$(\theta)$-stability & $89.02^o$  & 0 & $(2.3147,1.8796,1.5822)$ \\
			4  & Strong A$(\theta)$-stability & $87.34^o$ & 0.27 & $(3.9396,2.4506,2.2271,2.0612)$ \\
			\hline
		\end{tabular}
	\end{table}
	
	Given the $\omega_j$ parameters of Table \ref{tableTRK}, our goal is to provide a complete study of the stability of TASE-RK methods, without assuming that the matrix $A$ necessarily corresponds to the exact Jacobian of the problem.

	\section{Stability diagrams of TASE-RK methods}\label{sec3}
	In this section, we define the so-called stability diagrams related to the TASE-RK methods, proving some geometric properties that will be very useful for the results of the rest of the paper.
	From now on, we denote by ${\bf{0}}$ a column vector with all zeros of appropriate size.
	
	As already mentioned in Section \ref{sec1}, let us consider the test problem \eqref{test_ODEs}, performing the splitting of the Jacobian $J=A+B$. Since $J$ is known from the problem and $A$ is set by the user in the TASE-RK method \eqref{TASE_RK}, the matrix $B$ is consequently defined as $B:=J-A$. The IVP \eqref{test_ODEs} thus becomes
	\begin{equation}\label{EqTest1mat}
	\vec{u}_t=A \vec{u} + B \vec{u}+ \vec{g}(t),
	\end{equation}
	with $\vec{u}(t_0)=\vec{u}_0$. Given a problem to solve, our goal is to identify choices of $A$ guaranteeing the unconditional stability of TASE-RK methods, or at least stability for acceptable values of the step-size $k$. In the latter case, fixing $A$, we want to provide a technique for computing a value of $k$ which leads to the stability of the method.
	
	\begin{definition}\label{def1}{\normalfont{Conditional and unconditional stability of TASE-RK methods.}}
		A TASE-RK method \eqref{TASE_RK} is (conditionally) stable if, when applied to \eqref{EqTest1mat} with $\vec{g}={\bf{0}}$, there exist two constants $C\geq 0$ and $k^*>0$ such that
		\begin{equation*}
			||\vec{u}_n||\leq C, \quad \forall n=0,\ldots,N, \quad \forall k\leq k^*.
		\end{equation*}
		If $k^*=\infty$, the method is unconditionally stable.
		The constant $C$ could depend on the matrices $A$ and $B$ and the vector ${\vec{\omega}}$, but not on the time step-size and the index $n$.
	\end{definition}

	As done in \cite{Rosales20172336}, we consider the case where $A$ is real, Hermitian (symmetric in the real case) and negative definite, i.e. with all strictly negative eigenvalues:
	\begin{equation*}
		A^T=A, \qquad \langle \vec{u},A\vec{u} \rangle <0, \quad \forall \vec{u} \neq {\bf{0}}, \quad \vec{u} \in \mathbb{C}^d.
	\end{equation*}
	Now and in the following, we use the standard notation $\langle \vec{x},\vec{y} \rangle=\sum_{j=1}^d \bar{x}_j y_j$, $\vec{x},\vec{y} \in \mathbb{C}^d$, with $||\vec{x}||^2=\langle \vec{x},\vec{x} \rangle$. We underline that this hypothesis on $A$ is not restrictive. Indeed, once $A$ is fixed, with $J$ known from the problem, $B$ is defined accordingly as $B:=J-A$, and therefore the sum of $A$ and $B$ returns the Jacobian. Rather, choosing $A$ symmetric and negative definite is convenient, since in TASE-RK methods we have to solve linear systems with coefficient matrices of the form $(I-\omega_j k A)$, $\omega_j>0$, $k>0$, at each step.
	We recall that several existing direct or iterative methods for linear systems have great efficiency and robustness when the related coefficient matrices are symmetric and positive definite \cite[Ch. IV]{Trefethen}, as in this case.

\begin{figure}[h!]
		\centering
		\includegraphics[scale=0.435]{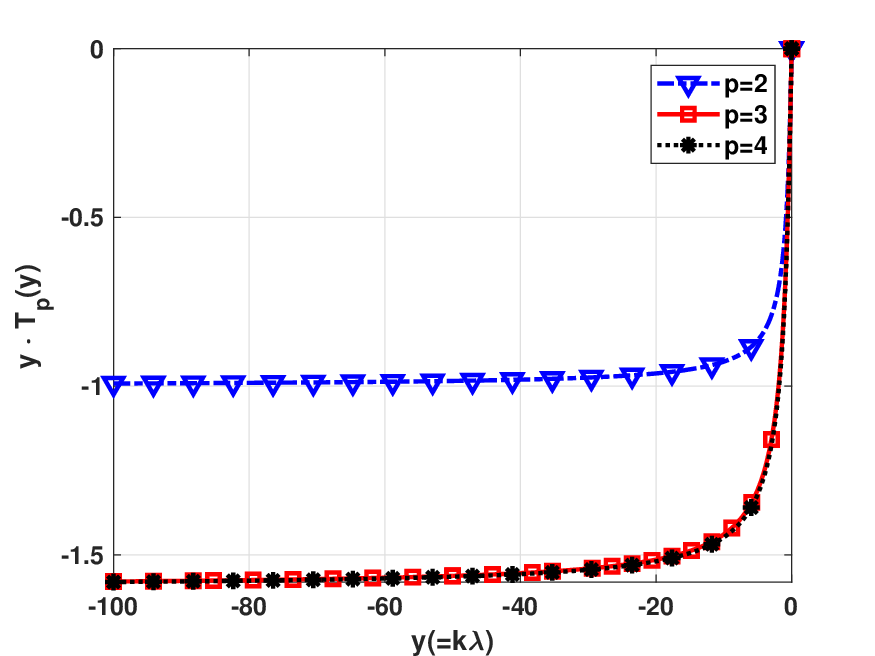}
		\caption{Function $\hat{T}_p(y)(=yT_p(y))$, for $p=2,3,4$. For each $p$, $\hat{T}_p$ is an increasing function whose values are all within the interval $[\hat{T}_p^*,0]$, where $\hat{T}_2^*\approx -1$, $\hat{T}_3^*\approx -1.5961$, $\hat{T}_4^*\approx -1.5961$.}\label{Fig11}
	\end{figure}

	Let $\eta:=\lambda+\gamma$ be a complex parameter such that Re$(\eta)\leq0$, representative of the eigenvalues of $ J $, where $\lambda \in \mathbb{R}^-$ and $\gamma \in \mathbb{C}$ play the role of the eigenvalues of $A$ and $B$, respectively. From \eqref{EqTest1mat}, we therefore write the scalar equation
	\begin{equation}\label{EqTest1}
	u_t=\lambda u + \gamma u + g(t).
	\end{equation}
	In this way we are implicitly assuming the simultaneous diagonalizability of $A$ and $B$; however, this procedure is necessary to define the stability diagrams of the TASE-RK methods, which will be useful also when we remove this hypothesis.
	Using \eqref{EqTest1} as test equation, the stability function \eqref{Stab_TASE} of $s$-stage TASE-RK methods with order $p=s\leq 4$ becomes
	\begin{equation}\label{Stab_TASE_IMEX}
	RT_p(k\lambda,k\gamma)=1+T_p(k\lambda) (k\lambda + k\gamma)+\frac{1}{2}(T_p(k\lambda) (k\lambda + k\gamma))^2+\ldots+\frac{1}{p!}(T_p(k\lambda)(k\lambda+k\gamma))^p.
	\end{equation}
	For simplicity of notation, from now on we avoid highlighting the dependency on ${\vec{\omega}}$ in the TASE operator, and therefore also in the stability function.
	Let us define $\mu:=\lambda^{-1} \gamma$, and make the change of variable $\gamma\rightarrow \mu \lambda$. Setting $y:=k\lambda\ (< 0)$, from \eqref{Stab_TASE_IMEX} we get the transformed stability function $\tilde{R}T_p(y,\mu):=RT_p(y,y \mu)$:
	\begin{equation}\label{Rtilde}
	\tilde{R}T_p(y,\mu)=1 + \hat{T}_p(y)(1+\mu) + \frac{1}{2} \hat{T}_p^2(y)(1+\mu)^2 + \ldots + \frac{1}{p!} \hat{T}_p^p(y)(1+\mu)^p.
	\end{equation}
	Here, we have denoted $\hat{T}_p(y):=yT_p(y)$, $y< 0$. In Figure \ref{Fig11}, we plot $\hat{T}_p$, for $p=2,3,4$. We will use the function $\hat{T}_p$ in the rest of the paper in several circumstances.
	
	\begin{definition}\label{def2}{\normalfont{Stability diagrams of TASE-RK methods.}}
		The (conditional) stability diagram of a TASE-RK method (with $\tilde{R}T_p$ as in \eqref{Rtilde}) is defined as
		\begin{equation*}\label{D}
			\mathcal{D}_{y,p}:=\{\mu \in \mathbb{C}: |\tilde{R}T_p(y,\mu)|\leq1,  \text{ with Re}(\mu)\geq-1 \}.
		\end{equation*}
		For $y\rightarrow -\infty$ ($k\rightarrow\infty$), we get the unconditional stability diagram, denoted by $\mathcal{D}_{\infty,p}$.
	\end{definition}
	The restriction Re$(\mu)\geq- 1$ comes from the assumption Re$(\eta)\leq0$ (recall that $\eta:=\lambda+\gamma$ and $\mu:=\lambda^{-1} \gamma$).
	Note that, in general, $\mathcal{D}_{y,p}$ depends on the eigenvalues of $A$ and $B$ and on the step-size $k$. However, $\mathcal{D}_{\infty,p}$ is independent of the step-size $k$. In the rest of this section, we derive some geometric properties of $\mathcal{D}_{y,p}$ that will be used to prove the conditional and unconditional stability properties of TASE-RK methods, whether $A$ and $B$ are simultaneously diagonalizable or not.

\begin{figure}[h!]
	\centering
	\includegraphics[scale=0.435]{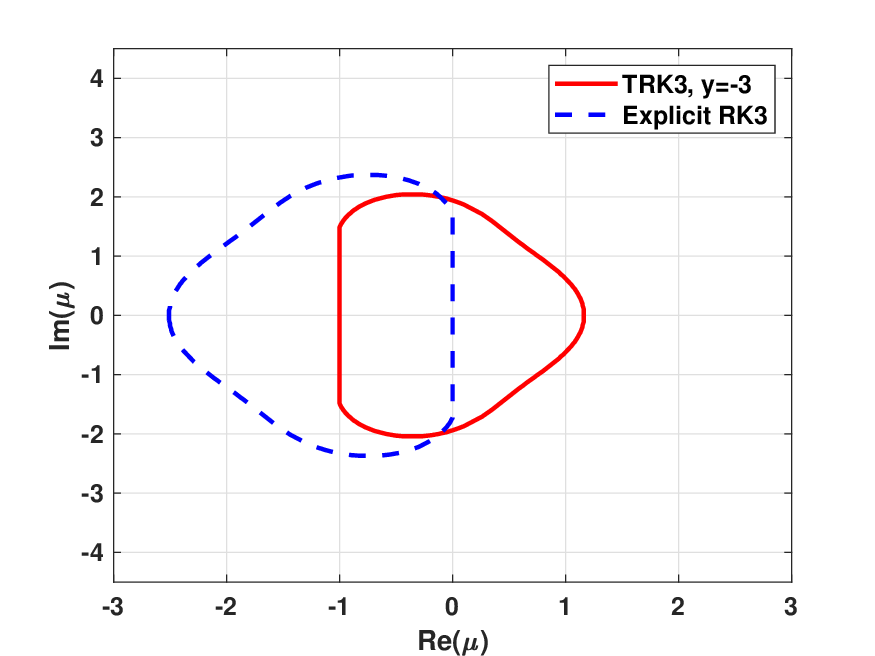} \includegraphics[scale=0.435]{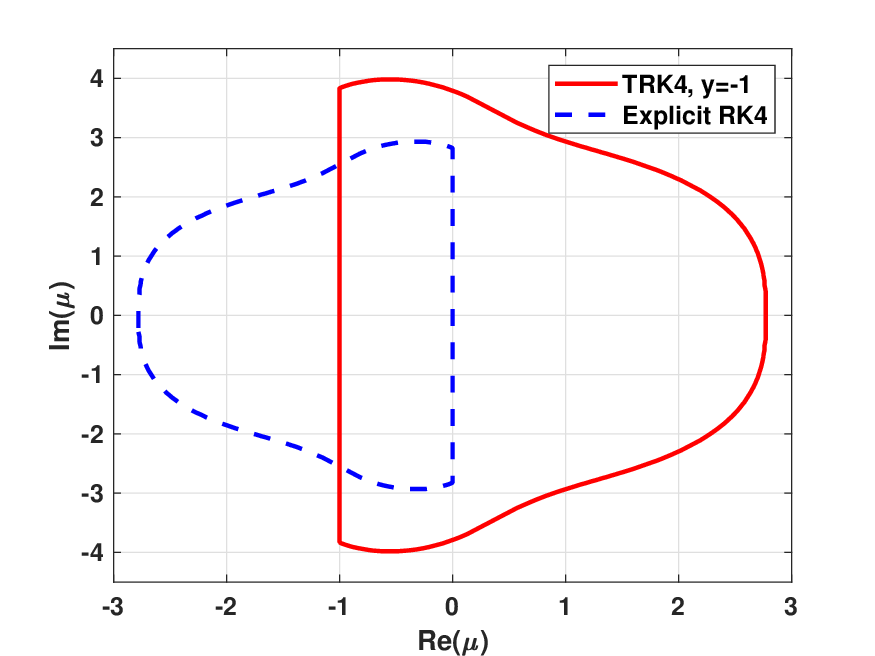}
	\caption{Stability region $\mathcal{R}_p$ (inside the dashed line) of explicit $p$-stage order-$p$ RK (RK$p$) methods versus stability diagram $\mathcal{ D }_{y,p}$ (inside the continuous line, fixing some values of $y$) of $p$-stage order-$p$ TASE-RK (TRK$p$) methods for $p=3$ (left) and $p=4$ (right).}\label{Fig2}
\end{figure}

	\begin{remark}\label{remark2}
		It is well known that the stability function $R_p(z)$ ($=RT_p(z)$ in \eqref{Stab_TASE} for $T_p(z)=1$) of an explicit $s$-stage RK method with order $p=s\leq 4$ is independent of its coefficients. The corresponding stability region, given by all $z \in \mathbb{C}$, Re$(z)\leq 0$, such that $|R_p(z)|\leq 1$, therefore depends only on $p$. Let us denote it with $\mathcal{R}_p$. For $p=2,3,4$, $\mathcal{R}_p$ has the following properties \cite[Subsec. 238]{Butcherbook}:
		\begin{itemize}
			\item $\mathcal{R}_p$ is simply connected (roughly speaking, there are no holes within $\mathcal{R}_p$);
			\item $0 \in \partial \mathcal{R}_p$, where in this context $0$ is the origin (in the complex plane) and the symbol $\partial$ stands for the boundary of the region;
			\item under the crucial assumption Re$(z)\leq 0$ (which means we do not consider the right-half complex plane), $\overline{0 P} \subset \mathcal{R}_p$, $\forall P \in \partial \mathcal{R}_p$, i.e. every segment whose endpoints are the origin and a generic $ P \in\partial \mathcal{R}_p$ is entirely contained in $\mathcal{R}_p$.
		\end{itemize}
		In Figure \ref{Fig2}, we report $\mathcal{R}_p$ for $p=3, 4$ (on the left and on the right, respectively; $\mathcal{R}_p$ is the region inside the dashed line).
	\end{remark}

	\begin{lemma}\label{lemma3.1}
		$\mathcal{D}_{y,p}$ is obtained from $\mathcal{R}_p$ through the composition of the following two elementary transformations:
		\begin{itemize}
			\item an homothety of coefficient $\hat{T}^{-1}_p(y)$;
			\item a shift on the real axis of ‘length’ $-1$.
		\end{itemize}
		Formally, this can be written as: $\mathcal{ D}_{y,p} = -1 + \hat{T}^{-1}_p(y) \mathcal{ R}_p$.
		\begin{proof}
			The proof is immediate by observing that, from \eqref{Rtilde}, setting $\tilde{\mu}=\mu+1$ and $\tilde{z}=\hat{T}_p(y)\tilde{\mu}$, we get
			\begin{equation*}
				\tilde{R}T_p(y,\mu)=R_p(\tilde{z})=1 + \tilde{z} + \frac{1}{2} \tilde{z}^2 + \ldots + \frac{1}{p!} \tilde{z}^p.
			\end{equation*}
			Thus, a point $\mu \in \mathcal {D} _{y,p} $ if and only if $\tilde{z} \in \mathcal{ R}_p$, where, for the calculations shown above, $\tilde{z}$ is obtained from $\mu$ through the composition of a shift on the real axis of ‘length’ $+1$ and an homothety of coefficient $\hat{T}_p(y)$. Working backwards, this concludes the proof.
		\end{proof}
	\end{lemma}
	
	Figure \ref {Fig2} shows a graphic example of the fact that, setting generic values of $y$, the region $\mathcal{D}_{y,p}$ is actually the transformation of $\mathcal{R}_p$ as described in Lemma \ref{lemma3.1}. Note that $\mathcal{D}_{y,p}$ is rotated $180^o$ relative to $\mathcal{R}_p$ (relative to the third dimension in 3D space), and this happens since the homothety coefficient $\hat{T}_p(y) \leq 0$, $\forall y < 0$ (see Figure \ref{Fig11}).
	
	\begin{corollary}\label{lemma3.2}
		For $p=2,3,4$, and any $y< 0$:
		\begin{itemize}
			\item $\mathcal{D}_{y,p}$ is simply connected;
			\item $P_{-1} \in \partial \mathcal{ D }_{y,p}$, where $P_{-1}$ is the point of coordinates $(-1,0)$ in the complex plane (the coordinates refer to the real and imaginary axis, respectively);
			\item $\overline{P_{-1} P} \subset \mathcal{ D }_{y,p}$, $\forall P \in \partial \mathcal{ D }_{y,p}$.
		\end{itemize}
		\begin{proof}
			The proof is straightforward from Remark \ref{remark2} and Lemma \ref{lemma3.1}.
		\end{proof}
	\end{corollary}
	Using the above lemma and corollary, we prove below the main result of this section.
	\begin{theorem}\label{prop2}
		For each $y_1<y_2< 0$ and $p=2,3,4$, it holds that $\mathcal{D}_{y_1,p} \subseteq \mathcal{D}_{y_2,p}$.
		\begin{proof}
			The proof is carried out by means of the following steps:
			\begin{itemize}
				\item we get the analytic expression of the points of the boundary of $\mathcal{D}_{y,p}$, solving for equality $\tilde{R}T_p(y,\mu)=e^{i\theta}$;
				\item we shift by coefficient $+1$ on the real axis the diagram $\mathcal{D}_{y,p}$ to simplify the calculations, denoting the region thus obtained by $\mathcal{\tilde{D}}_{y,p}$;
				\item fixing $y_1<y_2< 0$, we prove that $\mathcal{\tilde{D}}_{y_1,p} \subseteq \mathcal{\tilde{D}}_{y_2,p}$, and therefore $\mathcal{{D}}_{y_1,p} \subseteq \mathcal{{D}}_{y_2,p}$.
			\end{itemize}
			
			{ \raggedright \textbf{- Case $\bf p=2$.} }
			Solving for equality $\tilde{R}T_2(y,\mu)=e^{i\theta}$, we get
			\begin{equation*}
				\begin{split}
					& \mu_l=-1+\frac{c_l}{\hat{T}_2(y)}, \text{ with } c_l=-1-\sqrt{-1+2e^{i\theta}}, \text{ Re}(\mu_l)\geq -1, \\
					& \mu_r = -1+\frac{c_r}{\hat{T}_2(y)}, \text{ with } c_r=-1+\sqrt{-1+2e^{i\theta}}, \text{ Re}(\mu_r)\geq -1.
				\end{split}
			\end{equation*}
			
			Varying the angle $\theta$, $\mu_l$ and $\mu_r$ represent the points of $\partial \mathcal{ D }_{y,2}$. The points of $\partial \mathcal{\tilde{ D} }_{y,2}$ are therefore given by $\tilde{\mu}_j(y)=c_j/\hat{T}_2(y)$, $j=l,r$.
			
			Corollary \ref{lemma3.2} implies that $0 \in \partial \mathcal{ \tilde{D} }_{y,2}$ and $\overline{ 0 P} \subset \mathcal{ \tilde{D} }_{y,2}$, for any $P \in \partial \mathcal{\tilde{ D }}_{y,2}$. Note that $|c_j/\hat{T}_2(y_2)|\geq |c_j/\hat{T}_2(y_1)|$, $j=l,r$, since $\hat{T}_2(y)$ is a negative increasing function for each $y<0$ (see Figure \ref{Fig11}). Thus, given any $\theta$, we have $\overline{0 \tilde{\mu}_j(y_1)} \subseteq \overline{0 \tilde{\mu}_j(y_2)} (\subset \mathcal{ \tilde{D} }_{y_2,2})$.
			
			Since $\tilde{\mu}_j(y_1) \in \partial {\mathcal{ \tilde{D} }_{y_1,2}}$ (by definition), for the simple connectivity of ${\mathcal{\tilde{ D }}_{y,2}}$ (for each $y< 0$), we have that $\partial {\mathcal{ \tilde{D} }_{y_1,2}} \subset {\mathcal{ \tilde{D} }_{y_2,2}}$, and so $ {\mathcal{ \tilde{D} }_{y_1,2}} \subseteq {\mathcal{ \tilde{D} }_{y_2,2}}$.
			
			{ \raggedright  \textbf{- Case $\bf p=3$.}	}
			Solving for equality $\tilde{R}T_3(y,\mu)=e^{i\theta}$, we get
			\begin{equation*}
				\begin{split}
					&\mu_l=-1+\frac{c_l}{\hat{T}_3(y)}, \text{ with } c_l =\frac{1}{c}-c-1 , \text{ Re}(\mu_l)\geq -1, \\& \mu_{r_{\pm}}=-1+\frac{c_{r_{\pm}}}{\hat{T}_3(y)}, \text{ with } c_{r_{\pm}}=\frac{-1\pm i \sqrt{3}}{c}+c(1\pm i \sqrt{3})-2, \text{ Re}(\mu_{r_{\pm}})\geq -1,
				\end{split}
			\end{equation*}
			where $c=(1-3 e^{i\theta}+\sqrt{2-6ie^{i\theta}+9e^{2i\theta}})^{1/3}$. Repeating the same steps as the previous case, exploiting that the points of $\partial \mathcal{\tilde{D}}_{y,3}$ are all of the form $\mu_j=c_j/\hat{T}_3(y)$,
			where $\hat{T}_3(y)$, $y< 0$, is a negative increasing function (see Figure \ref{Fig11}), the theorem follows.
			
			{\raggedright  \textbf{- Case $\bf p=4$.}	}
			Solving for equality $\tilde{R}T_4(y,\mu)=e^{i\theta}$, we get
			\begin{equation*}
				\begin{split}
					& \mu_{l_{\pm}}=-1+\frac{c_{l_{\pm}}}{2\hat{T}_4(y)}, \text{ with } c_{l_{\pm}} =-2+2c_2\pm \sqrt{\varphi(c_1,c_2)} , \text{ Re}(\mu_{l_{\pm}})\geq -1, \\ & \mu_{r_{\pm}}=-1+\frac{c_{r_{\pm}}}{\hat{T}_4(y)}, \text{ with } c_{r_{\pm}}=-2-2c_2\pm\sqrt{\varphi(c_1,c_2)+\frac{16}{c_2}}, \text{ Re}(\mu_{r_{\pm}})\geq -1,
				\end{split}
			\end{equation*}
			where $c_1=-1 + 6 e^{i \theta} -
			\sqrt{-3 + 12 e^{i \theta} - 12 e^{2 i \theta} + 32 e^{3 i \theta}}$, $c_2=\sqrt{-1 + (c_1/2)^{1/3} +}$  $\overline{(2/c_1)^{1/3} (1 - 2 e^{i \theta})}$,
			and
			$\varphi(c_1,c_2)=-2^{5/3} c_1^{1/3} - 8 (1 + c_2)/c_2 +
			2^{7/3} (-1 + 2 e^{i \theta})/c_1^{1/3}$.
			Repeating the steps of the two previous cases, the theorem follows.
		\end{proof}
	\end{theorem}
	
\begin{figure}[h!]
		\centering
		\includegraphics[scale=0.435]{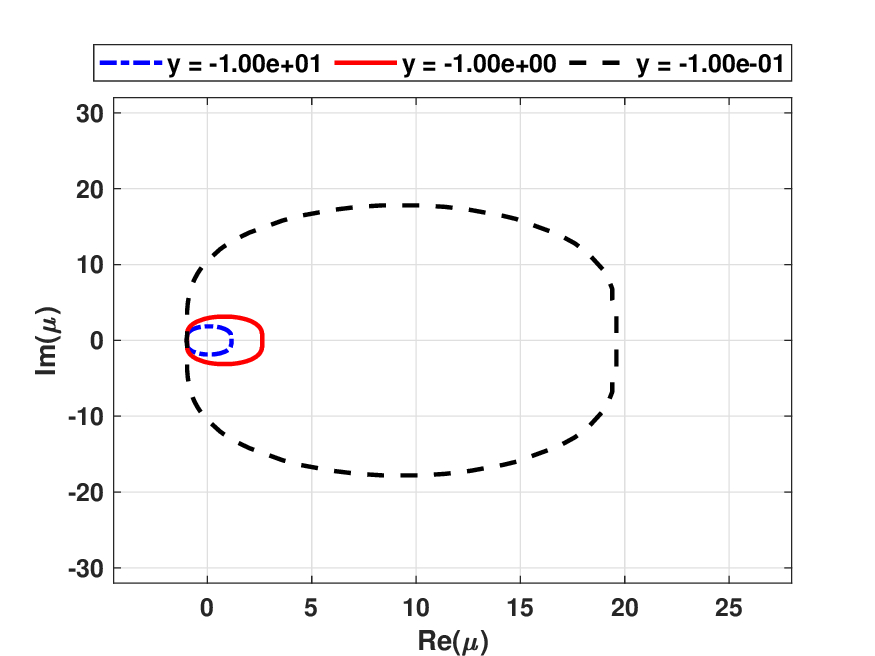} \includegraphics[scale=0.435]{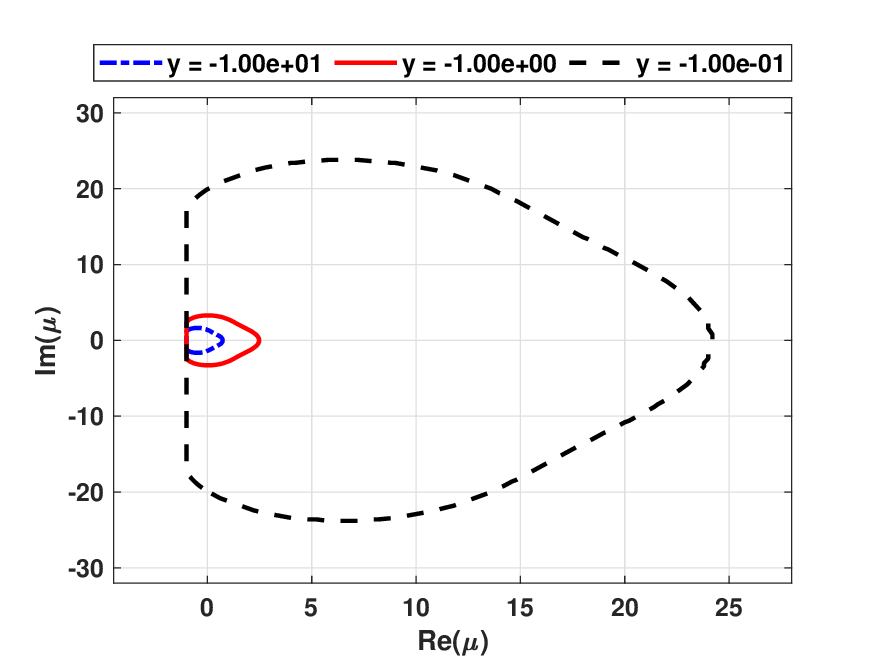}
		\includegraphics[scale=0.435]{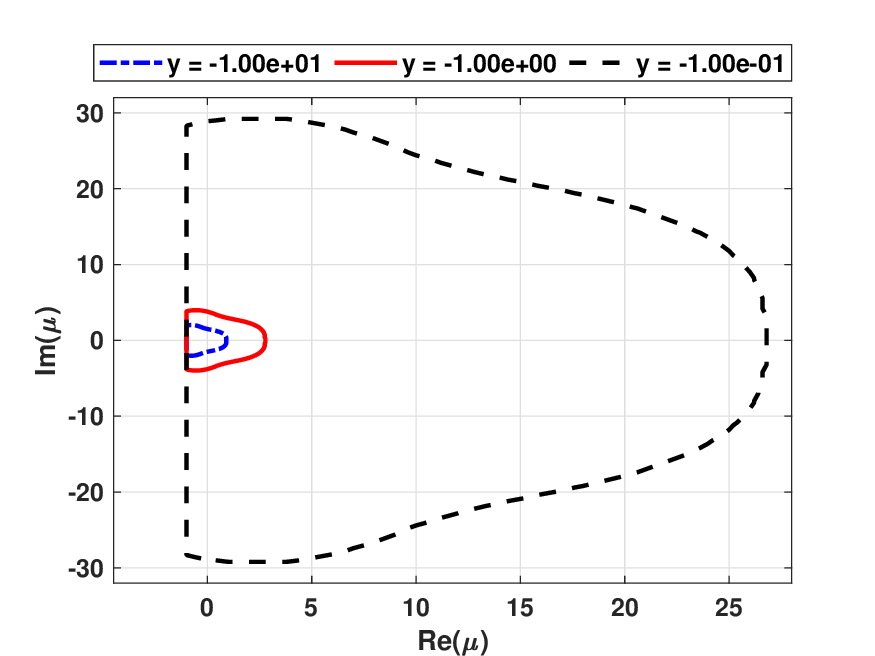} \includegraphics[scale=0.435]{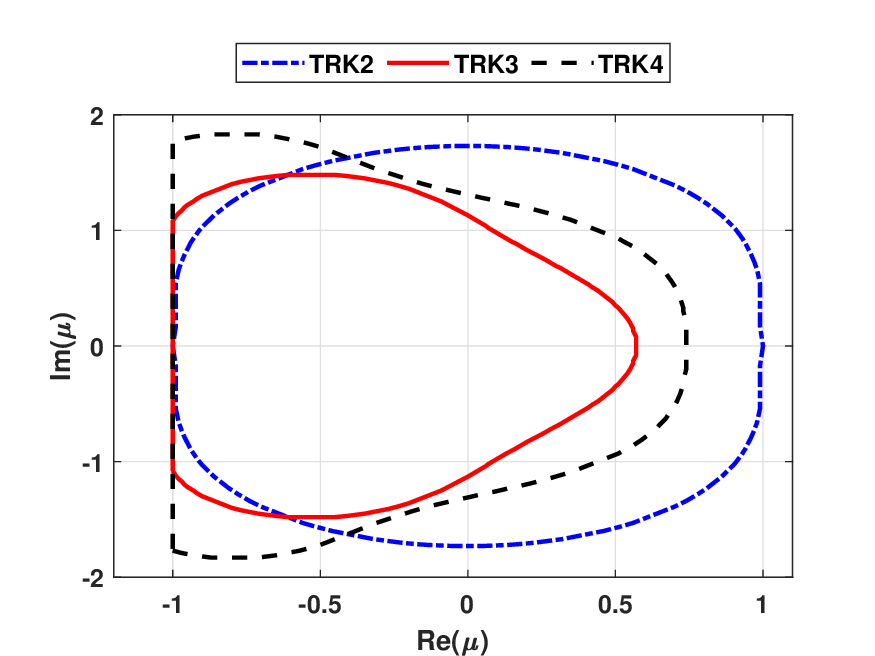}
		\caption{Stability diagrams $\mathcal{D}_{y,p}$ for $p=2,3,4$ (top left, top right, bottom left, respectively) for several values of $y$; $\mathcal{D}_{\infty,p}$ for $p=2,3,4$ (bottom, right).}\label{Fig3}
	\end{figure}

	We report in Figure \ref{Fig3} (top and bottom left) the stability diagrams $\mathcal{ D }_{y,p}$, $p=2,3,4$, in correspondence of several values of $ y $. In this way, it is possible to have a graphic view of what has been proved in Theorem \ref{prop2}.
	
	\begin{remark}\label{rem3}
		Theorem \ref{prop2} implies that $\mathcal{D}_{\infty,p}$ is the smallest stability diagram, i.e. $\mathcal{D}_{\infty,p} \subseteq \mathcal{D}_{y,p}\ \forall y< 0$. Figure \ref{Fig3} (bottom, right) shows $\mathcal{D}_{\infty,p}$ for $p=2,3,4$. Such diagrams will be useful to prove the unconditional stability properties of TASE-RK methods.
	\end{remark}

	\section{Stability of TASE-RK methods in case of simultaneous diagonalizability}\label{sec4}
	
	In this section we prove some necessary and/or sufficient conditions on the stability of TASE-RK methods assuming that $A$ and $B$ are simultaneously diagonalizable. In order not to make the statements of propositions and theorems redundant, we will avoid repeating the validity of this assumption every time. Let us first introduce the set of the so-called generalized eigenvalues of the splitting $(A,B)$:
	\begin{equation*}\label{gen_eig_set}
		\upmu(A,B):=\{\mu \in \mathbb{C} : \mu A \vec{v} = B \vec{v},\ \vec{v} \neq {\bf{0}}\}.
	\end{equation*}
	Note that the elements of $\upmu(A,B)$ are just the eigenvalues of $A^{-1}B$. From now on we use $\sigma(X)$ to denote the eigenvalues of a matrix $X$.
	
	The next proposition provides a necessary and sufficient condition for the stability of TASE-RK methods when the one-to-one correspondence between the eigenvalues of $A$ and $B$ is known.
	
	\begin{proposition}\label{prop1}
		Let $A \vec{v} = \lambda_i \vec{v}$, $B \vec{v} = \gamma_i \vec{v}$, with $i \in \{1,\ldots,d\}$, and $\vec{v} \in \mathbb{C}^d$. Furthermore, let $\mu_i= \lambda_i^{-1} \gamma_i$. Given a value of the time step-size $k>0$, a TASE-RK method with $s=p \leq 4$ is stable if and only if $\mu_i \in \mathcal{D}_{y_i,p}\ \forall i=1,\ldots,d$, where $y_i=k \lambda_i$.
	\end{proposition}	
	The proof of this result has been postponed to the appendix. 	
	
	\begin{remark}\label{rem2}
		{\textit{Computation of the step-size $k^*$ ensuring the stability of a TASE-RK method.}}
		We derive a procedure for the computation of $k^*$ in Definition \ref{def1} when $\mu \in \mathbb{R}$, for TASE-RK methods (with $p=s\leq 4$). To this end, we need to characterize the points of $\partial \mathcal{ D}_{y,p} $ on the real axis, given that the largest value of $k^*$ ensuring stability is obtained for $|\tilde{R}T_p(y,\mu)|=1$. Theorem \ref{prop2} helps us in this case, since in the related proof we have determined the expression of each point of $\partial \mathcal{ D}_{y,p} $ for $p=2,3,4$.
		
		{ \raggedright \textbf{- Case $\bf p=2$.} }
		The left and right points of $\partial \mathcal{ D}_{y,2} $ on the real axis are, respectively,
		\begin{equation*}
			\mu_l = -1, \qquad \mu_r(y)=-1-\frac{2}{\hat{T}_2(y)}.
		\end{equation*}
		Therefore, $\mu_l=-1$ is a fixed point of $ \partial \mathcal {D} _ {y,2} $, while $\mu_r(y) \geq 1\ \forall y < 0$, since $\hat{T}_2(y)$ is a negative increasing function bounded by $-1$ and $0$ (see Figure \ref{Fig11}). Therefore, naturally, if $\mu_i\leq 1\ \forall i=1,\ldots,d$, the TASE-RK method with $s=p=2$ is unconditionally stable (given that, for all $i$, $\mu_i \in \mathcal{D}_{y,2}\ \forall y< 0$ in this case, and therefore $\mu_i \in \mathcal{D}_{\infty,2}$). Otherwise, using Proposition \ref{prop1}, remembering that $y_i=k\lambda_i$ and $\mu_i= \lambda_i^{-1} \gamma_i$, by solving for $k$ the equations $\mu_i=\mu_r(y_i)$, $i=1,\ldots,d$, we get
		\begin{equation*}\label{kstarTRK2}
			k^* = \min_{i}\Bigg\{\frac{-8 + \mu_i - \sqrt{28 + 20 \mu_i + \mu_i^2}}{9 \lambda_i (-1 + \mu_i)}\Bigg\}, \qquad \mu_i> 1, \quad \lambda_i<0.
		\end{equation*}

		{ \raggedright \textbf{- Case $\bf p=3$.} } The left and right points of $\partial \mathcal{ D}_{y,3} $ on the real axis are, respectively,
		\begin{equation*}
			\mu_l = -1, \qquad \mu_r(y)=-1+\frac{ (-4 + \sqrt{17})^{1/3}-1-1/(-4 + \sqrt{17})^{1/3}}{\hat{T}_3(y)}.
		\end{equation*}
		Thus, $\mu_l=-1$ is a fixed point of $ \partial \mathcal {D} _ {y,3} $, while $\mu_r(y) \geq \mu_3^*$, with $\mu_3^*\approx 0.5743$, $\forall y < 0$, since $\hat{T}_3(y)$ is a negative increasing function bounded by $\hat{T}_3^*\approx -1.5961$ and $0$ (see Figure \ref{Fig11}). Similarly to the previous case, if $\mu_i\leq \mu_3^*\ \forall i=1,\ldots,d$, the TASE-RK method with $s=p=3$ is unconditionally stable. Otherwise, $k^*$ can be derived by solving for $k$ the equations $\mu_i=\mu_r(y_i)$, $i=1,\ldots,d$, then choosing the smallest $k$.

		{ \raggedright \textbf{- Case $\bf p=4$.} }
		The left and right points of $\partial \mathcal{ D}_{y,4} $ on the real axis are, respectively,
		\begin{equation*}
			\mu_l = -1, \qquad \mu_r(y)=-1+\frac{2^{2/3}(-43 + 9 \sqrt{29})^{1/3}-4 - 10 (2/(-43 + 9 \sqrt{29}))^{1/3}}{3\hat{T}_4(y)}.
		\end{equation*}
		Repeating the same steps as in the previous cases, we get that, if  $\mu_i \leq \mu_4^*$, with $\mu_4^*\approx 0.7445$, $\forall i=1,\ldots,d$, the TASE-RK method with $s=p=4$ is unconditionally stable. Otherwise, $k^*$ can be derived as already described above.
	\end{remark}

	Below is an example of application of Proposition \ref{prop1} and Remark \ref{rem2}.

	\begin{example}\label{test1}
		Let us consider the problem $\vec{u}_t=A\vec{u} + B\vec{u} + 10 \cdot  \mathbf{1}$, where from now on $\mathbf{1}$ represents the column vector of ones of appropriate length, with
		\begin{equation}\label{ex1m}
		A=\begin{pmatrix}[1.5]
		-40 & 30 & 30 \\
		30 & -\frac{71}{2} & -\frac{69}{2} \\
		30 & -\frac{69}{2} & -\frac{71}{2} \\
		\end{pmatrix}, \qquad B=\begin{pmatrix}[1.5]
		-\frac{74}{3} & \frac{38}{3} & \frac{38}{3} \\
		\frac{38}{3} & -\frac{233}{12} & -\frac{215}{12}\\
		\frac{38}{3} & -\frac{215}{12} & -\frac{233}{12} \\
		\end{pmatrix}.
		\end{equation}
		The matrices $A$ and $B$ are simultaneously diagonalizable. The corresponding eigenvalues are $(\lambda_1,\lambda_2,\lambda_3)=(-100,-10,-1)$, $(\gamma_1,\gamma_2,\gamma_3)=(-50,-12,-\frac{3}{2})$, respectively. Knowing the one-to-one correspondence between the eigenvalues of $A$ and $B$, we can compute $\mu_i=\lambda_i^{-1} \gamma_i$: $(\mu_1,\mu_2,\mu_3)=(\frac{1}{2},\frac{6}{5},\frac{3}{2})$. To have TASE-RK stability (see Proposition \ref{prop1}), the step-size $k$ must be chosen in such a way that $\mu_i \in \mathcal{D}_{y_i}$, $y_i=k\lambda_i, \ \forall i$. By means of Remark \ref{rem2}, we can compute $k^*$, i.e. the largest value of $k$ that guarantees stability. Denoting by TRK$p$ the $s$-stage TASE-RK method with order $p=s$, we have e.g. $k^*\approx 7.8390e-01$ for TRK2, and $k^*\approx 2.8428e-01$ for TRK3.

		\begin{figure}[h!]
			\centering
			\includegraphics[scale=0.435]{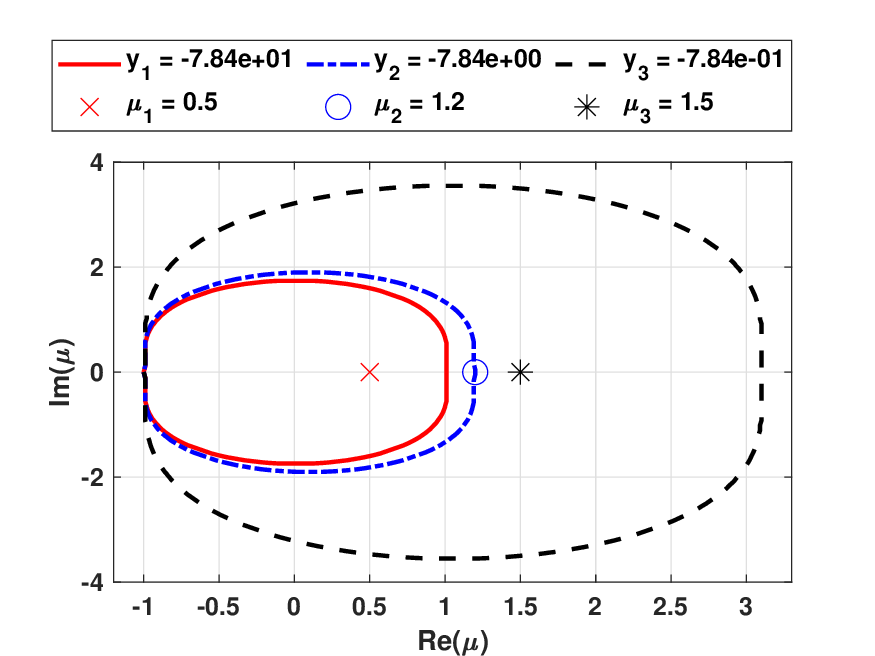} \includegraphics[scale=0.435]{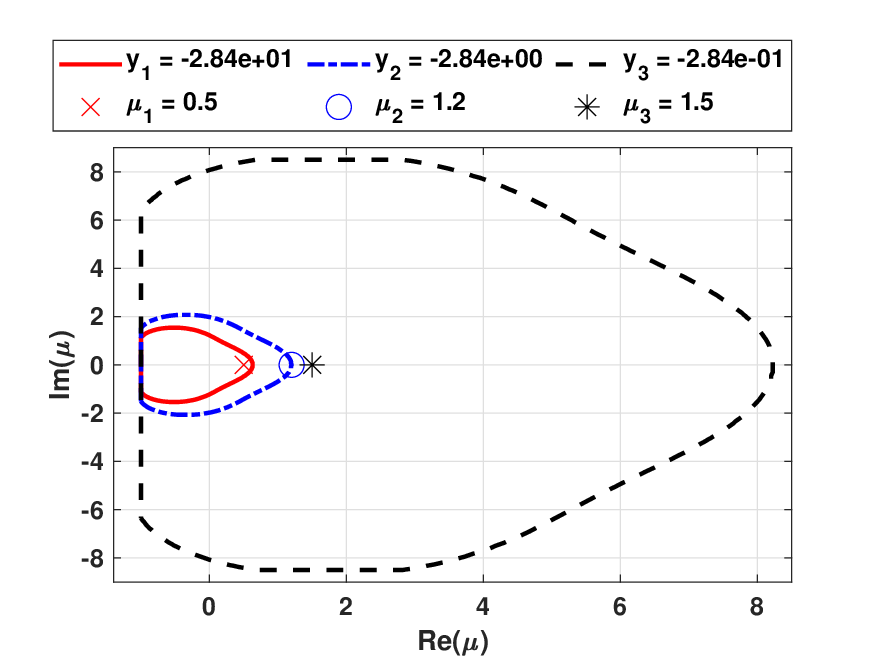}
			\caption{Stability diagrams $\mathcal{D}_{y_i,p}$ of TRK2 (left) and TRK3 (right) methods, for several values of $y_i=k\lambda_i$, i.e. $(\lambda_1,\lambda_2,\lambda_3)=(-100,-10,-1)$, $k=7.8390e-01$ for TRK2, $k=2.8428e-01$ for TRK3; the points plotted inside the figures correspond to $(\mu_1,\mu_2,\mu_3)=(\frac{1}{2},\frac{6}{5},\frac{3}{2})$.}\label{Fig1Test1}
		\end{figure}
		
		Figure \ref{Fig1Test1} shows the $\mathcal{D}_{y_i,p}$ regions and the corresponding $\mu_i$, $i=1,2,3$, for TRK2 and TRK3, given these values of $k$ ($k=k^*$). Note that for both methods $\mu_1$ and $\mu_3$ are inside the related stability diagrams, while $\mu_2$ is just on the boundary of $\mathcal{D}_{y_2,p}$. Thus, according to Proposition \ref{prop1}, TRK2 and TRK3 are stable for the mentioned step values. Furthermore, by choosing $k>k^*$, the methods become unstable. This is confirmed in Table \ref{Ex1}, where we report the errors due to the use of the TRK2 and TRK3, choosing initial condition $\vec{u}(t_0)=(2,3,1) \cdot 10^2$, with $t_0=0$. We applied the TRK2 and TRK3 methods also with exact Jacobian $J=A+B$ in the TASE operator; in this case, in accordance with the results of Table \ref{tableTRK}, there are no stability issues.
		
		\begin{table}[h!]
			\centering
			\caption{Relative errors in norm two by TRK2 and TRK3 at the point $t_e=30$, related to numerical solution of $\vec{u}_t=A\vec{u} + B\vec{u} + 10 \cdot \mathbf{1}$ (fixing the coefficient matrix in the TASE operator corresponding to both the exact Jacobian $J$ and $A$), with $A$ and $B$ as in \eqref{ex1m} and initial condition $\vec{u}(0)=(2,3,1) \cdot 10^2$.} \label{Ex1}
			%{\color{red}
			\begin{tabular}{c || c c | c c }
				\hline
				$k$ & TRK2 ($J$) & TRK2 ($A$) & TRK3 ($J$) & TRK3 ($A$) \\
				\hline
				$1.8750e+00$ & $8.1916e-03$ & $2.6260e+03$ & $3.2074e-10$ & $1.1479e+10$ \\
				$9.3750e-01$ & $3.4523e-07$ & $1.1609e+03$ & $2.3285e-15$ & $5.3503e+14$ \\
				$4.6875e-01$ & $2.1246e-13$  & $2.5721e-01$ & $3.0794e-16$  & $1.3881e+16$ \\
				$2.3438e-01$ & $6.1587e-16$ & $1.5684e-15$ & $6.5109e-16$ & $9.5785e-13$ \\
				\hline
			\end{tabular}
		\end{table}

	\end{example}

	Proposition \ref{prop1} and Remark \ref{rem2} can be very useful, but they presuppose knowledge of the one-to-one correspondence between the eigenvalues of $A$ and $B$. For this reason, we derive below a sufficient condition for the stability of TASE-RK methods that does not require such knowledge.

	\begin{theorem}\label{prop3}
		Let $\sigma(A)=\{\lambda_i, \ i=1,\ldots,d: \lambda_1 \leq \lambda_2 \leq \ldots \leq \lambda_d<0 \}$. Given a value of the time step-size $k>0$, a TASE-RK method with $s=p \leq 4$ is stable if $\upmu(A,B)\subseteq \mathcal{ D }_{y_1,p} $, where $y_1=k\lambda_1$.
		\begin{proof}
			The proof follows from Proposition \ref{prop1} and Theorem \ref{prop2}. Indeed, Proposition \ref{prop1} states that the condition $\mu_i \in \mathcal{D}_{y_i,p}$ ($\forall  i$) guarantees the stability of the TASE-RK method. Since $\mu_i= \lambda_i^{-1} \gamma_i$ are just the generalized eigenvalues of $(A,B)$, this means that stability is guaranteed if $\upmu(A,B)\subseteq \mathcal{ D }_{y_i,p}\ \forall i=1,\ldots,d$. From Theorem \ref{prop2}, $\mathcal{D}_{y_1,p} \subseteq \mathcal{D}_{y_i,p}\ \forall i=1,\ldots,d$ ($y_i=k\lambda_i$), and therefore the proof follows.
		\end{proof}
	\end{theorem}

	\begin{remark}

		Theorem \ref{prop3} provides a more restrictive condition than Proposition \ref{prop1} on the choice of the step-size $k$ for TASE-RK stability.
		Let us look, for instance, at Test \ref{test1}. In that case, we know that the largest step-sizes we can choose to have TRK2 and TRK3 stability are $k^*=7.8390e-01$, $k^*=2.8428e-01$, respectively. Indeed, for both methods, $\mu_i(=\lambda_i^{-1} \gamma_i) \in \mathcal{ D }_{y_i,p}$, $y_i=k\lambda_i, \ \forall i$. Actually, $\mu_2$ lies just on the boundary of $\mathcal{D}_{y_2,p}$, see Figure \ref{Fig1Test1}. Therefore, given the eigenvalues $(\lambda_1,\lambda_2,\lambda_3)=(-100,-10,-1)$ of $A$, the one that leads to the step-size restriction is $\lambda_2$ (and not $\lambda_1$, which is the smallest). If we had not known the correspondence between each $\lambda_i$ and $\mu_i$, we would have used Theorem \ref{prop3} instead of Proposition \ref{prop1}, which leads to the stability requirement $\mu_i \in \mathcal{D}_{y_1 ,p} \ \forall i=1,2,3$. This implies choosing $k=k^*/10$ in this case (since $\lambda_1=10\lambda_2$) to ensure the stability of the TASE-RK methods.
	\end{remark}    	
	Although it provides more restrictions on $k$ for TASE-RK stability, using Theorem \ref{prop3} is much simpler than Proposition \ref{prop1}.
	
	To conclude this section, we provide a necessary and sufficient condition for the unconditional stability of TASE-RK methods.
	
	\begin{theorem}\label{th42}
		A TASE-RK method with $s=p \leq 4$ is unconditionally stable if and only if $\upmu(A,B)\subseteq \mathcal{D}_{\infty,p} $.
		\begin{proof}
			Let us first prove that requiring $\upmu(A,B)\subseteq \mathcal{D}_{\infty,p} $  is sufficient for the unconditional stability. From Theorem \ref{prop2} (see also Remark \ref{rem3}), $\mathcal{D}_{\infty,p} \subseteq \mathcal{D}_{y,p}$, for each $y< 0$. Exploiting Theorem \ref{prop3}, we therefore get the proof.
			
			Instead, let us now assume that the TASE-RK method is unconditionally stable. Using Definition \ref{def1}, this means that
			\begin{equation*}
				||\vec{u}_n||\leq C, \quad \forall n=0,\ldots,N, \quad \forall k>0,
			\end{equation*}
			where $\vec{u}_n$ is the time-marching solution of the test problem \eqref{EqTest1mat} with $\vec{g}={\bf{0}}$, obtained through the TASE-RK method \eqref{TASE_RK}. By means of the same steps as the proof of Proposition \ref{prop1} up to recurrence \eqref{recurrence}, we therefore have that
			\begin{equation*}
				|\tilde{R}T_p(\infty,\mu_i)|\leq 1,   \quad \forall i=1,\ldots,d, \quad \mu_i=\lambda_i^{-1} \gamma_i \in \upmu(A,B).
			\end{equation*}
			Thus, by Definition \ref{def2}, $\mu_i \in \mathcal{D}_{\infty,p}\ \forall i=1,\ldots,d$. This concludes the proof.
		\end{proof}
	\end{theorem}
	
	In this section, we have analyzed the conditional and unconditional stability of the TASE-RK methods when $A$ and $B$ are simultaneously diagonalizable,  also showing a simple example.
	In the next section, we study the stability of TASE-RK methods by removing the assumption of simultaneous diagonalizability of $A$ and $B$.

	\section{Stability of TASE-RK methods in the general case}\label{sec5}
	To analyze the stability of TASE-RK methods in general cases, we exploit some connections between $\mathcal{D}_{y,p}$ and the Field Of Values (FOV) of a matrix $X$ \cite[Ch. I]{HornJohnsonFOV}, defined as
	\begin{equation}\label{FOVdef}
	\mathcal{W}(X):=\{\langle \vec{x}, X \vec{x} \rangle : ||\vec{x}||^2 = 1, \ \vec{x} \in \mathbb{C}^d \}.
	\end{equation}
	In particular, we consider the set
	\begin{equation}\label{FOVexpr}
	\mathcal{W}_q(-A,B):=\{\langle \vec{v}, (-A)^{q-1}B\vec{v} \rangle : \langle \vec{v}, (-A)^{q}\vec{v} \rangle = 1, \ \vec{v} \in \mathbb{C}^d, \ q\in \mathbb{R}\}.
	\end{equation}
	With $\vec{v}=(-A)^{-q/2}\vec{x}$, it can be written as
$
		\mathcal{W}_q(-A,B)= \mathcal{W} ((-A)^{q/2-1} B (-A)^{-q/2})$,
	where $\mathcal{W}$ is exactly the FOV defined in \eqref{FOVdef}. Note that we are considering $-A$ (which is positive definite) to make the definition \eqref{FOVexpr} well-posed.
	
	The connections between the stability diagrams and the FOV have already been exploited in \cite{Rosales20172336} for the study of the unconditional stability of a class of LMMs. In this section, we take inspiration from this study to analyze the conditional and unconditional stability of TASE-RK methods.
	
	Let us recall the analytic and geometric properties of the FOV that we will use.
	\begin{remark}\label{remFOV}
		Let $X_1$, $X_2$, be two square matrices of size $d$, and $w_1$, $w_2 \ \in \mathbb{C}$:
		\begin{enumerate}
			\item $\mathcal{ W }( X_1) + \mathcal{ W }( X_2) = \{x_1+x_2: x_1 \in \mathcal{ W }(X_1),\ x_2 \in \mathcal{W}(X_2)\} $;
			\item $\mathcal{ W }( X_1) \mathcal{ W }( X_2) = \{x_1x_2: x_1 \in \mathcal{ W }(X_1),\ x_2 \in \mathcal{W}(X_2)\} $;
			\item $\mathcal{ W }(w_1 X_1+w_2 I)=w_1 \mathcal{ W }(X_1)+w_2 $ \cite[p. 9, Properties 1.2.3 and 1.2.4]{HornJohnsonFOV};
			\item $\sigma(X_1) \subset \mathcal{ W }(X_1)$ \cite[p. 10, Property 1.2.6]{HornJohnsonFOV};
			\item if $X_1$ is normal, then  $\mathcal{ W }(X_1) = \text{co}(\sigma(X_1))$ \cite[p. 11, Property 1.2.9]{HornJohnsonFOV};  co$(\sigma(X))$ denotes the convex hull of the eigenvalues of the matrix $X$;
			\item as a direct consequence of \cite[Corollary 1.7.7]{HornJohnsonFOV}, if $X_1$ is positive semi-definite, then $\sigma(X_1 X_2) \subseteq \mathcal{W}(X_1) \mathcal{ W }(X_2)$.
		\end{enumerate}
			There exist several codes to efficiently compute the FOV; we will use the MATLAB \texttt{chebfun} routine \cite{chebfun}, based on a classical algorithm proposed by Johnson \cite{JohnsonSINUM1978}.
	\end{remark}

	Note that, since the assumption of simultaneous diagonalizability of $A$ and $B$ has been removed, we can no longer simplify the calculations by making the transition from the test IVP \eqref{EqTest1mat} to the scalar equation \eqref{EqTest1} for the study of the stability. Let us use Definition \ref{def1}, thus applying the TASE-RK method \eqref{TASE_RK} to the problem \eqref{EqTest1mat} with $\vec{g}={\bf{0}}$. This leads to the recurrence
	\begin{equation}\label{recurrence_matrix}
	\vec{u}_{n+1}= (RT_p(kA,kB))^{n+1} \vec{u}_{n},
	\end{equation}
	where $RT_p$ is here the natural extension of \eqref{Stab_TASE_IMEX} to the multidimensional case. Denoting by $\rho(X)$ the spectral radius of a matrix $X$, from \eqref{recurrence_matrix} we therefore conclude that the stability condition of TASE-RK methods in the general case reads
	\begin{equation}\label{stabcond_TASE_nocommute}
	\rho(RT_p(kA,kB))\leq 1.
	\end{equation}
	
	The following lemma will be useful to establish two subsequent theorems about the stability of TASE-RK methods.
	\begin{lemma}\label{lemmaC}
		Let $ \sigma(X)=\{c_i, \ i=1,\ldots,d \}$, with $X$ square matrix. Given some coefficients $w_q, \ q=0,\ldots,p$, it holds that
		\begin{equation*}\label{thesis_lemmaC}
			\sigma(\sum_{q=0}^pw_q X^q) = \{\sum_{q=0}^pw_q c_i^q, \ i=1,\ldots,d \}.
		\end{equation*}
	\end{lemma}
	The proof of this result has been postponed to the appendix.
	
	We are now ready to prove the next result, which links the stability of TASE-RK methods to the region $\mathcal{R}_p$ (recall that with $\mathcal{R}_p$ we indicate the stability region of an explicit $s$-stage RK method with order $p=s\leq 4$). This result will be employed to prove some necessary or sufficient conditions on the stability of TASE-RK methods.
	\begin{theorem}\label{th40}
		Given a value of the time step-size $k > 0$, a TASE-RK method with $s=p \leq 4$ is stable if and only if $\sigma(\mathcal{Z}) \subseteq \mathcal{ R}_p$, where $\mathcal{Z}:=kT_p(kA)\ (A+B)$.
		\begin{proof}
			Let $ \sigma(\mathcal{Z})=\{\zeta_i, \ i=1,\ldots,d \}$. If a TASE-RK method is stable, then
			\begin{equation}\label{rel1}
			\bigg|1+\zeta_i+\frac{1}{2}\zeta_i^2+\ldots + \frac{1}{p!}\zeta_i^p \bigg| \leq 1, \quad \forall i=1,\ldots,d.
			\end{equation}
			Through Lemma \ref{lemmaC}, taking $w_q=1/q!$, $q=0,\ldots,p$, note that \eqref{rel1} is equivalent to
			\begin{equation}\label{rel2}
			\bigg|\sigma\bigg(I+\mathcal{Z}+\frac{1}{2}\mathcal{Z}^2+\ldots + \frac{1}{p!}\mathcal{Z}^p\bigg) \bigg| \leq 1.
			\end{equation}
			Observe that \eqref{rel2} corresponds to the stability condition \eqref{stabcond_TASE_nocommute} of TASE-RK methods.
		\end{proof}
	\end{theorem}
	
	We prove below a sufficient condition for the stability of TASE-RK methods, and then two conditions (one sufficient, one necessary) for unconditional stability. We will use the following simple lemma during the proofs.
		
	\begin{lemma}\label{lemma5.4}
		Let $ \sigma(A) = \{\lambda_i,\ i=1,\ldots,d\}$. It holds that
		$$
		\sigma \big( (kA) \ T_p(kA) \big)=\{\hat{T}_p(y_i) , \ i=1,\ldots,d\},$$
		where we recall the notation $y_i=k\lambda_i$ and $\hat{T}_p(y):=yT_p(y)\ \forall y< 0$.
		\begin{proof}
			Observe that $ \sigma(T_p(kA)) = \{T_p(y_i),\ i=1,\ldots,d\}$, since $T_p(kA)$ is the sum of commuting matrices of the form $\beta_j(I-\omega_jkA)^{-1}$, whose eigenvalues are given by $\beta_j/(1-\omega_j y_i)$, $i=1,\ldots,d$. Then, exploiting the simultaneous diagonalizability of $kA$ and $T_p(kA)$, the theorem follows.
		\end{proof}
	\end{lemma}
	
	\begin{theorem}\label{cond2}
		Let $\sigma(A)=\{\lambda_i, \ i=1,\ldots,d:\lambda_1 \leq \lambda_2 \leq \ldots \leq \lambda_d<0\} $. Given a value of the time step-size $k>0$, a TASE-RK method with $s=p \leq 4$ is stable if $\exists\ q \in \mathbb{R}: \mathcal{W}_q(-A,B) \subseteq -\mathcal{D}_{y_1,p}$, where $y_1=k\lambda_1$.
		\begin{proof}
			From Theorem \ref{th40}, we need to prove $\sigma(\mathcal{Z})\subseteq \mathcal{R}_p$, where $\mathcal{Z}:=kT_p(kA) \ (A+B)$. Let $\zeta \in \sigma(\mathcal{ Z })$, therefore $\exists\ \vec{v}\neq {\bf 0}: \mathcal{ Z }\vec{v}=\zeta \vec{v}$. Given $q \in \mathbb{R}$, we can write
			\begin{equation*}
				\mathcal{Z}(-A)^{-q/2}(-A)^{q/2}\vec{v}=(-A)^{-q/2}(-A)^{q/2}\zeta \vec{v}.
			\end{equation*}
			Multiplying both sides from the left by $(-A)^{q/2}$, then setting $\vec{x}=(-A)^{q/2}\vec{v}$, we have $\zeta \in \sigma((-A)^{q/2}\mathcal{Z}(-A)^{-q/2})$; we will therefore prove that $\sigma((-A)^{q/2}\mathcal{Z}(-A)^{-q/2}) \subseteq \mathcal{ R }_p$.
			
			Let us decompose $\mathcal{Z}$ into the product of two matrices,
			as follows:
			\begin{equation*}
				\mathcal{Z}=M N, \qquad M=kT_p(kA)\ (-A), \quad N=-I+(-A)^{-1}B.
			\end{equation*}
			Observing that the matrices $(-A)^{q/2}$ and $M$ commute, we have that
			\begin{equation}\label{rel11}
			\sigma((-A)^{q/2}\mathcal{Z}(-A)^{-q/2})=\sigma(M(-A)^{q/2}N(-A)^{-q/2}).
			\end{equation}
			Since $M$ is symmetric positive definite,  by means of Remark \ref{remFOV}, item 6, we get
			\begin{equation}\label{rel21}
			\sigma(M(-A)^{q/2}N(-A)^{-q/2}) \subseteq \mathcal{W}(M)\mathcal{W}((-A)^{q/2}N(-A)^{-q/2}).
			\end{equation}
			Note that $\mathcal{W}(M)=\text{co}(\sigma(M))$, as $M$ is a normal matrix (see Remark \ref{remFOV}, item 4). From Lemma \ref{lemma5.4}, $\sigma(M)=\{-\hat{T}_p(y_i),\ i=1,\ldots,d\}$, $y_i=k\lambda_i$; co($\sigma(M)$) therefore corresponds (in the complex plane) to the segment $\overline{P_1 P_2}$, where $P_1=(-\hat{T}_p(y_d),0)$, $P_2=(-\hat{T}_p(y_1),0)$, as $\hat{T}_p(y)\leq 0 \ \forall y< 0$ (see Figure \ref{Fig11}). Since $ \overline{P_1 P_2}\subseteq \overline{ 0 P_2}$, we write
			\begin{equation}\label{rel31}
			\mathcal{W}(M)\mathcal{W}((-A)^{q/2}N(-A)^{-q/2}) \subseteq -\hat{T}_p(y_1) \mathcal{W}((-A)^{q/2}N(-A)^{-q/2}).
			\end{equation}
			
			Summarizing, from \eqref{rel11}, through \eqref{rel21} and \eqref{rel31}, so far we have obtained
			\begin{equation}\label{rel41}
			\sigma((-A)^{q/2}\mathcal{Z}(-A)^{-q/2}) \subseteq -\hat{T}_p(y_1) \mathcal{W}((-A)^{q/2}N(-A)^{-q/2}).
			\end{equation}
			Observing that
			\begin{equation*}
				(-A)^{q/2}N(-A)^{-q/2}=-I+(-A)^{q/2-1}B(-A)^{-q/2},
			\end{equation*}
			for the properties of the FOV (see Remark \ref{remFOV}, item 3) it holds
			\begin{equation*}
				\mathcal{W}((-A)^{q/2}N(-A)^{-q/2})=-1+\mathcal{ W }((-A)^{q/2-1}B(-A)^{-q/2})=-1+\mathcal{ W }_q(-A,B).
			\end{equation*}
			Since, by hypothesis, $\mathcal{ W }_q(-A,B) \subseteq -\mathcal{ D }_{y_1,p}$, from \eqref{rel41} we have
			\begin{equation*}
				\sigma((-A)^{q/2}\mathcal{Z}(-A)^{-q/2}) \subseteq \hat{T}_p(y_1)(1 +\mathcal{ D }_{y_1,p} ).
			\end{equation*}
			Noting that Lemma \ref{lemma3.1} implies $\mathcal{ R}_{p}=\hat{T}_p(y_1)(1 +\mathcal{ D }_{y_1,p} )$, the theorem follows.
		\end{proof}
	\end{theorem}
	
	Observe that the above sufficient condition for stability is the analogue of Theorem \ref{prop3} in the case in which $A$ and $B$ are not necessarily simultaneously diagonalizable; indeed, the only difference lies in the use of the FOV in place of the generalized eigenvalues. Obviously, the condition obtained in Theorem \ref{cond2} is more restrictive than the one given in Theorem \ref{prop3}, since $\upmu(-A,B)(=-\upmu(A,B))\subseteq \mathcal{ W }_q(-A,B) \ \forall q$ (from Remark \ref{remFOV}, item 4, using that $\sigma((-A)^{-1}B)=\sigma((-A)^{q/2-1}B(-A)^{-q/2})$).
	
	\begin{example}\label{test2}
		Let us consider the problem $\vec{u}_t=A\vec{u} + B\vec{u} + 10 \cdot  \mathbf{1}$, with
		\begin{equation}\label{ex2m}
		A=\begin{pmatrix}
		-10 & 0 & 0 \\
		0 & -4 & 0 \\
		0 & 0 & -30 \\
		\end{pmatrix}, \qquad B=\begin{pmatrix}
		-3 & 15 & 0 \\
		-15 & -3 & 0\\
		0 & 0 & -15 \\
		\end{pmatrix}.
		\end{equation}
		The matrices $A$ and $B$ are not simultaneously diagonalizable. Thus, to determine a condition that guarantees the stability of TASE-RK methods, we can use Theorem \ref{cond2}, choosing $k$ in such a way that $\mathcal{W}_q(-A,B)\subseteq -\mathcal{ D }_{y_1,p}$ $(y_1=k\lambda_1)$ for some $q$. In our case, see \eqref{ex2m}, $\lambda_1=-30$ (remember that, by $\lambda_i$, we denote the eigenvalues of $A$).
		
		However, observing that $A$ and $B$ are simultaneously block-diagonalizable, it is probably possible to require a less restrictive condition (although the paper does not present a rigorous proof of this statement; it will be the subject of future research). Indeed, $\lambda_1$ is associated with $\gamma_1=-15$, and so $\mu_1=\frac{1}{2}$. Note that $\mu_1=\frac{1}{2}$ belongs to the unconditional stability diagrams of all TASE-RK methods, see Figure \ref{Fig3} (bottom, right). Therefore, we apply Theorem \ref{cond2} only to the first $2 \times 2$ blocks of $A$ and $B$, consequently fixing $\lambda_1=-10$; since $\mathcal{D}_{-30k,p} \subseteq \mathcal{D}_{-10k,p}$ (see Theorem \ref{prop2}), this allows to obtain a less restrictive condition on the choice of $k$ (than using $\lambda_1=-30$). In Figure \ref{Figtest2}, we report $\mathcal{W}_1(-A,B)$ (left and right), $-\mathcal{D}_{-2.1,2}$ (left), $-\mathcal{D}_{-1.45,3}$ (right). Taking $\lambda_1=-10$, according to Theorem \ref{cond2}, we can therefore choose $k=2.1e-01$ for TRK2, and $k=1.45e-01$ for TRK3, being sure not to have stability issues. Solving the problem numerically, with the same initial conditions used in Test \ref{test1}, we have observed that the TRK2 and TRK3 methods are stable for $k<1$.
		\begin{figure}[h!]
			\centering
			\includegraphics[scale=0.435]{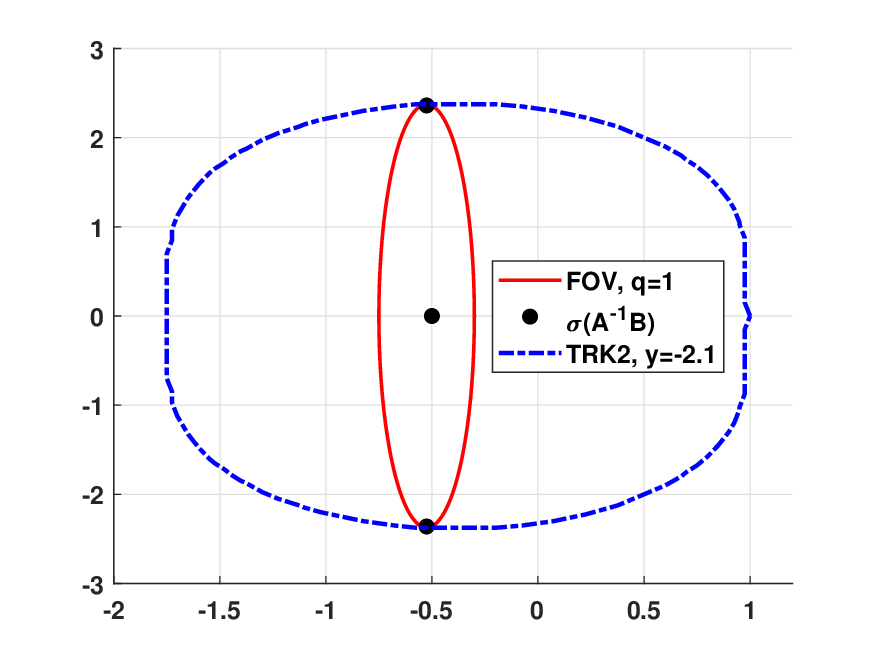} \includegraphics[scale=0.435]{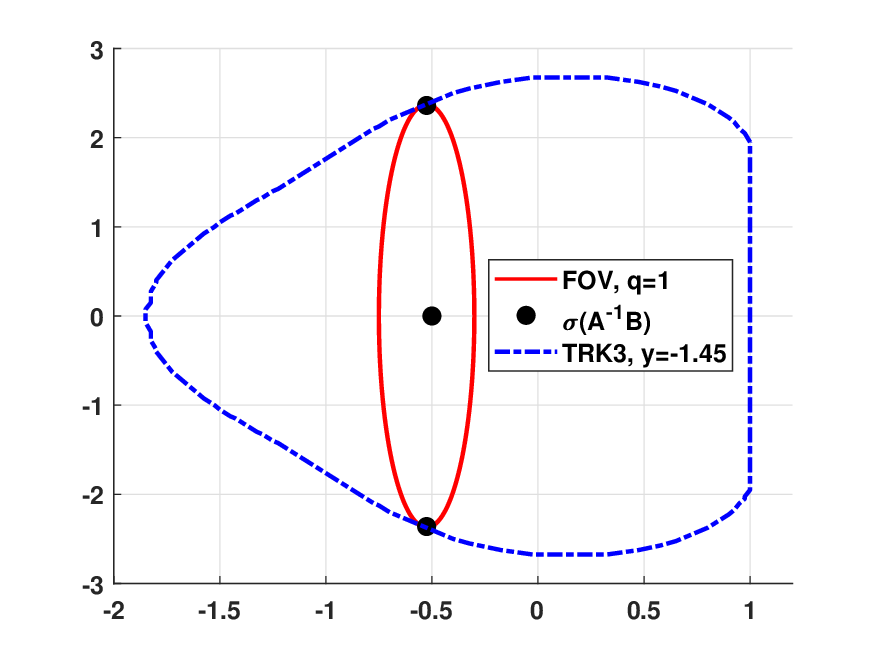}
			\caption{$\mathcal{W}_{1}(-A,B)$ (left and right), and $-\mathcal{D}_{y,2}$, $y=-2.1$ (left), $-\mathcal{D}_{y,3}$, $y=-1.45$ (right), with $A$ and $B$ as in \eqref{ex2m}.}\label{Figtest2}
		\end{figure}
	\end{example}

	Finally, we derive two conditions, one sufficient, one necessary, for the unconditional stability of TASE-RK methods.
	\begin{theorem}\label{th5.4}
		A TASE-RK method with $s=p \leq 4$ is unconditionally stable if $\exists\ q \in \mathbb{R}: \mathcal{W}_q(-A,B) \subseteq -\mathcal{D}_{\infty,p}$. If a TASE-RK method with $s=p \leq 4$ is unconditionally stable, then $\upmu(A,B) \subseteq \mathcal{D}_{\infty,p}$.
		\begin{proof}
			The sufficient condition naturally follows from Theorem \ref{cond2}.
			
			Thus, let us focus on the necessary condition. By hypothesis, setting $\mathcal{Z}:=kT_p(kA)\ (A+B)$, from \eqref{stabcond_TASE_nocommute} we have
			\begin{equation*}
				\rho	\bigg(I+\mathcal{Z}+\frac{1}{2}\mathcal{Z}^2+\ldots+\frac{1}{p!}\mathcal{Z}^p\bigg)\leq 1, \quad \forall k>0.
			\end{equation*}
			Using Lemma \ref{lemmaC}, this is equivalent to
			\begin{equation}\label{eqci}
			\begin{split}
			& \bigg| 1+c_i(k)+\frac{1}{2}c_i^2(k)+\ldots+\frac{1}{p!}c_i^p(k) \bigg|\leq 1, \quad \forall i=1,\ldots,d, \ \forall k>0, \\
			& \text{with } \ c_i(k)\in \sigma\big(kT_p(kA)\ (A+B)\big).
			\end{split}
			\end{equation}
			
			Writing $k$ as $kAA^{-1}$, then observing that $A^{-1}$ and $T_p(kA)$ commute, we have $kT_p(kA)\ (A+B)=\hat{T}_p(kA) \ (I+A^{-1}B)$. Here, $\hat{T}_p$ is the natural extension of the function $\hat{T}_p(y):=yT_p(y)$ to the multidimensional case.		
			Furthermore, from Lemma \ref{lemma5.4}, there exist two matrices $P$ and $D$ such that $\hat{T}_p(kA)=PDP^{-1}$, with $D=\text{diag}(\hat{T}_p(y_i))$, $y_i=k\lambda_i$, where $\lambda_i$ are, as usual, the eigenvalues of $A$. We therefore have that the eigenvalues $c_i$ in \eqref{eqci} satisfy
			\begin{equation}\label{cik}
			c_i(k)\in \sigma\big(P \ \text{diag}(\hat{T}_p(k \lambda_i)) \ P^{-1} \ (I+A^{-1}B)\big).
			\end{equation}
			Condition \eqref{eqci} holds for any $k>0$. Therefore, it remains valid even passing to the limit $k \rightarrow \infty$. Denoting with $\hat{T}_p^*=\lim_{k\rightarrow \infty} \hat{T}_p(k\lambda) \in \mathbb{R}$, from \eqref{cik}, condition \eqref{eqci} implies
			\begin{equation}\label{eqcstari}
			\begin{split}
			& \bigg| 1+{c_i}^*+\frac{1}{2}{c_i^*}^2+\ldots+\frac{1}{p!}{c_i^*}^p \bigg|\leq 1, \quad \forall i=1,\ldots,d,  \\
			& \text{with } \ {c_i}^* \in \sigma\big(\hat{T}_p^* \ (I+A^{-1}B) \big), \ \hat{T}_p^* \in \mathbb{R}.
			\end{split}
			\end{equation}
			Note that such a limit of $\hat{T}_p$ exists, and we know its value for $p=2,3,4$, see Figure \ref{Fig11}.
			
			Observing that $\sigma(\hat{T}_p^* \ (I+A^{-1}B))=\hat{T}_p^* \ (1+\upmu(A,B))$, the theorem follows, see Definition \ref{def2}. Indeed, \eqref{eqcstari} corresponds to  $\lim_{k \rightarrow \infty} |\tilde{R}T_p(k\lambda,\mu)|\leq 1$, with $\tilde{R}T_p$ as in \eqref{Rtilde} and $\mu$ any value in $\upmu(A,B)$.
		\end{proof}
	\end{theorem}
	Theorem \ref{th5.4} is the analogue of Theorem \ref{th42} in the previous section, and provides a very useful condition in practice. Indeed, if there exists a generalized eigenvalue of the splitting $(A,B)$ that does not belong to $\mathcal{D}_{\infty,p}$, we can already conclude that the TASE-RK method cannot be unconditionally stable.
	
	\begin{remark}\label{rem5}
		Of course, unconditional stability is a very strong demand and generally difficult to achieve. However, it may happen that, by multiplying the matrix $A$ by a suitably chosen positive parameter $\kappa$, thus consistently defining $B:=J-\kappa A$, it becomes possible for the TASE-RK method to be unconditionally stable. If we are already aware of $ \upmu (A, J) $ and/or $ \mathcal{W}_q (-A,J) $, we can directly derive $ \upmu (\kappa A, B) $ and/or $ \mathcal{W}_q (-\kappa A, B) $, exploiting the following relationships \cite[Rem. 5]{Seibold2019}:
		\begin{equation*}
			\upmu(\kappa A,B)=-1+\kappa^{-1}\upmu(A,J), \qquad \mathcal{W}_q(-\kappa A,B)=1+\kappa^{-1}\mathcal{ W}_q(-A,J).
		\end{equation*}
	\end{remark}

	\section{An application of the theory}\label{sec6}
	In this section we show an approach for applying the theory developed so far to problems in which the vector field $\vec{f}$ is characterized by a natural splitting into a linear and nonlinear part. For simplicity, let us first consider a scalar IVP of the type
	\begin{equation}\label{nonlinODE}
	u_t=\lambda u + g(u), \quad u(0)=u_0,
	\qquad \text{with} \quad \lambda \in \mathbb{R}^-,  \quad  g:\mathbb{R} \rightarrow \mathbb{R}, \quad t \in [0,t_e].
	\end{equation}
	We here assume that the function $g$ depends on $u$ and is nonlinear. Furthermore, we hypothesize that $g$ is continuous and differentiable.
	
	The Lagrange theorem (also known as mean value theorem) states that
	there exists $\xi=u(t^*)$, $t^*>0$, such that $g(u)-g(u_0)=g'(\xi)(u-u_0)$. Thus, the equation \eqref{nonlinODE} admits the rewriting
	\begin{equation}\label{mod_test}
	u_t= (\lambda+g'(\xi))u + k_g, \qquad \text{with} \quad k_g:=k_g(u_0,\xi) = g(u_0)-g'(\xi)u_0, \quad \xi \in \mathbb{R}.
	\end{equation}
	
	We can therefore use the theory proposed in the previous sections: the function $g'(\xi)$ now plays the role of the $\gamma$ parameter introduced in \eqref{EqTest1}. However, $\xi$ is unknown and difficult to determine. In the following example, we show one way to determine a value $\xi^*$ that provides a “safe bound”, in the context of the stability analysis, of $\xi$.
	
	\begin{example}\label{test3}
		Let us consider the scalar equation $u_t=\lambda u + g(u)$, $\lambda \in \mathbb{R}^-$, with $g(u)=\epsilon u(1-u)\in \mathbb{R}$, $u_0=u(0)=0$, which, according to \eqref{mod_test}, can be rewritten as
		\begin{equation} \label{FHNscalartest}
		u_t=(\lambda  + g'(\xi)) u, \qquad g'(\xi)=\epsilon(1-2\xi).
		\end{equation}
		
		Taking \eqref{FHNscalartest} as test equation, the stability condition of TASE-RK methods reads $|\tilde{R}T_p(y,\mu(\xi))| \leq 1$, $y=k\lambda$, $\mu(\xi)=\lambda^{-1} g'(\xi)$, with $\tilde{R}T_p$ as in \eqref{Rtilde}. By substituting the expressions for $y$ and $\mu(\xi)$ into $\tilde{R}T_p$, we obtain
		\begin{equation*}
			\tilde{R}T_p(z) =1 + z + \frac{1}{2} z^2 + \ldots + \frac{1}{p!} z^p, \qquad z=z(k,\lambda,\xi)=kT_p(k\lambda)(\lambda+g'(\xi)) \in \mathbb{R}.
		\end{equation*}
		Note that $\tilde{R}T_p$ takes the expression of the stability function of explicit $s$-stage RK methods of order $p=s\leq 4$. Requiring $|\tilde{R}T_p(z)| \leq 1$ leads to \cite[p. 101]{Butcherbook}
		\begin{equation*}\label{cp}
			-c_p \leq z \leq 0, \qquad \text{where}\quad c_2=2, \quad c_3\approx 2.5127, \quad c_4\approx 2.7853.
		\end{equation*}
		The condition $z \geq -c_p$, in this case, reads
		$
		g'(\xi) \geq -\frac{c_p}{kT_p(k\lambda)} - \lambda$.
		Being $\upsilon_{\text{LB}}, \ \upsilon_{\text{UB}}\  \in \mathbb{R}$ lower and upper bound for $u$, respectively, with $g'$ defined in \eqref{FHNscalartest}, we have $g'(\upsilon_{\text{UB}})\leq g'(\xi)\leq g'(\upsilon_{\text{LB}})$.
		Thus, requiring $k$ to satisfy
		\begin{equation*}
			g'(\upsilon_{\text{UB}}) \geq -\frac{c_p}{kT_p(k\lambda)} - \lambda =:\chi(p, k, \lambda),
		\end{equation*}
		leads to the stability of the method. That is, if the conditions $ |u_0 | \leq C, \ |\xi| \leq C \ \forall \xi $, imply $g'(\xi) \geq \chi(p, k, \lambda)$, then $ |u_1 | \leq C $, and therefore by induction $ |u_n | \leq C, \,  \forall n$.
		To conclude, one way to select “safe” value of $\xi^*$ in the context of  stability analysis can therefore be:
		$
		\xi^*=\arg \min g'(\xi)$.
	\end{example}

	In the non-scalar case, we have to consider IVPs of the form
	\begin{equation}\label{nonlinIVP}
	\vec{u}_t=A \vec{u} + \vec{g}(\vec{u}), \quad \vec{u}(0)=\vec{u}_0,
	\qquad \text{with} \quad \vec{g}:\mathbb{R}^d \rightarrow \mathbb{R}^d, \quad t \in [0,t_e].
	\end{equation}
	Writing the above nonlinear function as $\vec{g}(\vec{u})=(g_1(\vec{u}),g_2(\vec{u}),\ldots,g_d(\vec{u}))^T$, with $g_i(\vec{u})\in \mathbb{R}$, $i=1,\ldots,d$, we can exploit an extension of the Lagrange theorem, which states the existence of a vector $\vec{\xi}$ such that $g_i(\vec{u})-g_i(\vec{u}_0)=\langle \vec{\nabla} g_i(\vec{\xi}), \vec{u}-\vec{u}_0 \rangle$.
	Here, the operator $\vec{\nabla}$ is used to indicate the gradient of the functions $g_i$ at $\vec{\xi}$. Although complicated, it is therefore possible to make considerations similar to the scalar case.
	
	Many semi-discretized PDEs in space lead to IVPs \eqref{nonlinIVP}, with decoupled equations in $\vec{g}(\vec{u})$, i.e. such that $\vec{g}(\vec{u})=(g_1(u_1),g_2(u_2),\ldots,g_d(u_d))^T$, $u_i,\ g_i \ \in \mathbb{R}$, $i=1,\ldots,d$. In these cases, the study of stability simplifies, as explained through an example in the following subsection.

	\subsection{Fisher-Kolmogorov equation}\label{subsecFK}
	Let us consider the Fisher-Kolmogorov (FK) PDE \cite[Subsec. 13.2]{Murraybook}
	\begin{equation}\label{FKeq}
	u_t= D u_{xx} + \epsilon u(1 - u), \qquad D, \ \epsilon >0.
	\end{equation}
	The FK is a famous example of a nonlinear reaction-diffusion equation. It is related to the modeling of chemical reactions in which kinetics (nonlinear term) and diffusion are coupled.
	To reproduce the wave solutions of FK equation \eqref{FKeq}, we set $D=2 \cdot 10^{-2}, \epsilon=10^{-2}$, considering the following initial and boundary conditions:
	\begin{equation}\label{FKcond}
	\begin{split}
	& u_0:=u(x,t_0)=1+\frac{1}{2} e^{-x} \sin(\pi x), \qquad t_0=0, \quad x \in [-1,1],\\
	& u(-1,t)=u(1,t)=1.
	\end{split}
	\end{equation}
	We first apply the method of lines for the numerical solution of FK. Fixing the grid $\{x_m=-1+mh;\ m=0,\ldots,M;\ x_M=1\}$, then semi-discretizing the spatial derivative by means of finite differences of order two, we thus get from \eqref{FKeq} the system of ODEs
	\begin{equation*}
		\vec{u}_t=  A\vec{u} + \epsilon \vec{g}(\vec{u}) + \vec{u}_{\text{BC}},
	\end{equation*}
	where $\vec{u}_{\text{BC}}=(\frac{D}{h^2},0,\ldots,0,\frac{D}{h^2})^\text{T} \in \mathbb{R}^{M-1}$, and
	\begin{equation}\label{A}
	A= \frac{D}{h^2} \ \text{tridiag}(1,-2,1) \in \mathbb{R}^{M-1}, 	\quad \vec{g}(\vec{u})= (g_1(u_1),g_2(u_2),\ldots,g_d(u_d))^\text{T} \in \mathbb{R}^{M-1},
	\end{equation}
	with $g_i(u_i)=u_i(1-u_i) \in \mathbb{R}$. Applying Lagrange theorem, there exists a vector $\vec{\xi}=(\xi_1,\xi_2,\ldots,\xi_{M-1})$ such that the semi-discretized FK equation admits the rewriting
	\begin{equation}\label{B}
	\vec{u}_t=  (A+B(\vec{\xi}))\vec{u} + \vec{k}_{\vec{g}}, \qquad \text{with} \quad B(\vec{\xi})=\epsilon \ \text{diag}(1-2\xi_i)_{i=1}^{M-1}.
	\end{equation}
	Here, $\vec{k}_{\vec{g}}$ depends on $\vec{g}, \ \vec{u}_0, \ \vec{\xi}$ (see e.g. \eqref{mod_test}, which arises in the scalar context), and on the boundary conditions; let us assume its boundedness. % of $\vec{k}_{\vec{g}}$.
	
	Noting that the matrices $A$ and $B=B(\vec{\xi})$ do not commute, the stability of TASE-RK methods can be studied through Theorem \ref{cond2} (conditional stability) or Theorem \ref{th5.4} (unconditional stability). We therefore need to compute the FOV $\mathcal{W}_q(-A,B)$, setting a value of $q$, e.g. $q=1$, then verifying the condition on $k$ guaranteeing $\mathcal{W}_1(-A,B) \subseteq -\mathcal{D}_{y_1,p}$ ($y_1=k \lambda_1$), or if it happens that $\mathcal{W}_1(-A,B) \subseteq -\mathcal{D}_{\infty,p}$, respectively. The following proposition is very useful to this end.
	\begin{proposition}\label{prop6.1}
		With $A$ and $B=B(\vec{\xi})$ as in \eqref{A}, \eqref{B}, respectively, the FOV $\mathcal{W}_1(-A,B)$ takes only real values and satisfies
		\begin{equation}\label{proof6.1}
		\frac{\epsilon}{\ell} (1-2 \upsilon_{\text{UB}}) \leq \mathcal{ W}_1(-A,B) \leq \frac{\epsilon}{\ell} (1-2 \upsilon_{\text{LB}}),
		\end{equation}
		where $\ell \in [\frac{2D}{h^2}(1 + \cos(\frac{(M-1)\pi}{M})), \frac{2D}{h^2}(1 + \cos(\frac{\pi}{M})) ]$, $\upsilon_{\text{LB}}\leq \xi_i \leq \upsilon_{\text{UB}} \ \forall i=1,\dots,M-1$.
	\end{proposition}
	The proof of this result has been postponed to the appendix.
	
	From Theorem \ref{cond2}, to get the stability of the TASE-RK methods on the semi-discretized FK PDE, we have to choose $k$ such that $\mathcal{W}_1(-A,B)\subseteq -\mathcal{D}_{y_1,p}$, $y_1=k\lambda_1$, with $\lambda_1$ corresponding to smallest eigenvalue of the matrix $A$ in \eqref{A}. We can impose this condition easily, given that Proposition \ref{prop6.1} asserts that the FOV $\mathcal{W}_1$ takes only values in $\mathbb{R}$ (and provides us with bounds for them), and we know the analytical expression of the points of $\partial (-\mathcal{D}_{y_1,p})$ on the real axis from Remark \ref{rem2}. For each $p$, the expression of these points is given by $\mu=1-c_p/\hat{T}_p(y_1)$, $\mu=1$, respectively, with $y_1=k\lambda_1, \ c_p<0, \ \hat{T}_p(y_1)<0$. Therefore, from \eqref{proof6.1}, if the following holds, the TRK methods are stable:
	\begin{equation}\label{bounds}
	1 - \frac{c_p}{\hat{T}_p(k\lambda_1)}\leq \frac{\epsilon}{\ell} (1-2 \upsilon_{\text{UB}}), \qquad \frac{\epsilon}{\ell} (1-2 \upsilon_{\text{LB}}) \leq 1.
	\end{equation}
	
\begin{figure}[h!]
		\centering
		\includegraphics[scale=0.435]{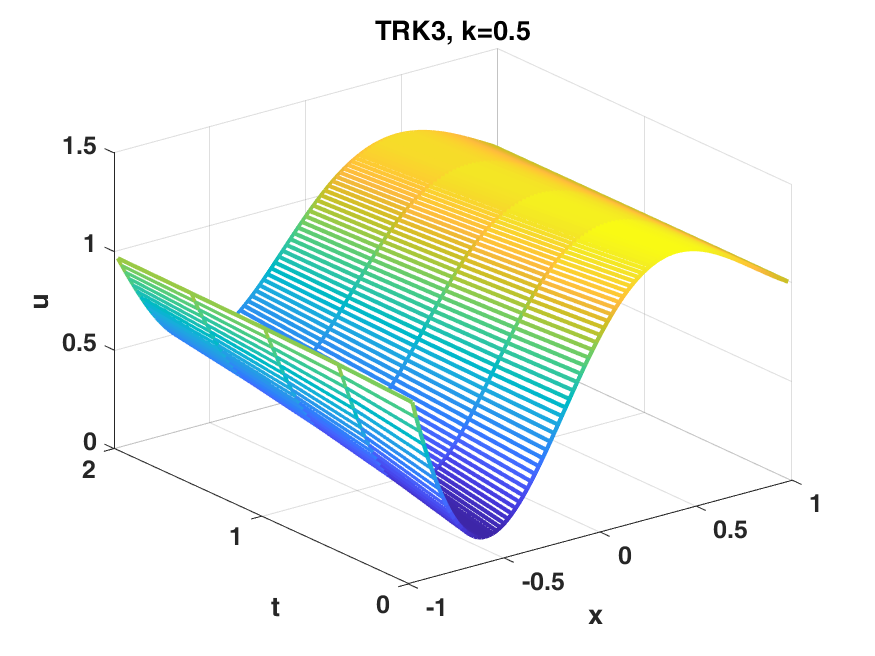}
		\caption{Numerical solution of the semi-discretized FK equation by the TRK3 method, with $k=0.5$, setting parameters, initial and boundary conditions as explained in Subsection \ref{subsecFK}.}\label{FKtest}
	\end{figure}

	Taking $\upsilon_{\text{LB}}=0, \ \upsilon_{\text{UB}}=\frac{3}{2}$ (see next Remark \ref{rem4}), and $M=100$, we get unconditional stability for all TASE-RK methods; see e.g. Figure \ref{FKtest}, where we report the numerical solution of the semi-discretized FK equation with initial and boundary conditions \eqref{FKcond}, by the TRK3 selecting a large value of $k$.

	\begin{remark}\label{rem4}
		In \eqref{bounds}, we have considered $\upsilon_{\text{LB}}=0$, and $\upsilon_{\text{UB}}=\frac{3}{2}$, which provide lower and upper bounds through $\arg \max g'(\xi)$, $\arg \min g'(\xi)$, at $(t_0, u_0)$, see \eqref{FKcond}, respectively. Semi-discretizing the FK PDE, we are therefore requiring $ \rho( RT_p ) <  1$ at $(t_0, \vec{u}_0)$ (thus, there exists a norm $|| . ||$, such that $ || RT_p || <  1$ at the first grid point). However, this does not guarantee that $ || \vec{u}_n ||_{\infty}  \leq  || \vec{u}_{n-1} ||_{\infty} $, nor $ || \vec{u}_n ||_{\infty}  \leq  \frac{3}{2} $. Hence, it is not possible to combine the maximum principle with induction principle to prove stability globally, unless we also check, after every iteration, that indeed $ || \vec{u}_n ||_{\infty}  \leq  \frac{3}{2} $. To prove bounds on the solution, Lyapunov theory could be used in some particular examples, but this idea has not yet been fully developed.
		
		Another way to make an estimate of  $\upsilon_{\text{LB}}$, $\upsilon_{\text{UB}}$, can e.g. be to perform convergence tests using some kind of manufactured solution, as explained in \cite[Subsec. 6.3]{Seibold2019}. In the numerical experiments, in order not to further burden the analysis of nonlinear stability and the procedure proposed to carry it out, we will use the information we have on the solution at $(t_0,\vec{u_0})$ (as done for FK); we will show that this still leads to satisfactory results. More accurate estimates of the mentioned bounds will be investigated in future work.
	\end{remark}

	\section{Numerical experiments}\label{sec7}
	This section is devoted to numerical experiments. The objective is twofold: on the one hand, we want to show the correctness of the theoretical results obtained in the manuscript; on the other hand, we want to show that thanks to these results the TASE-RK methods become very competitive and efficient. This last aspect will be analyzed by comparison with the following well-known and (some of them) recently introduced Rosenbrock (ROS) methods which, like TASE-RK, are a class of linearly implicit numerical schemes.
	\begin{itemize}
		\item \textbf{ROS2}: one-stage order-two ROS method derived in \cite[Eq. (1.26), p. 153]{Hundsdorferbook};
		\item \textbf{ROS3}: three-stage order-three ROS method derived \cite[Eq. (38), p. 573]{GonzalezPintoAMC2016};
		\item \textbf{ROS4}: four-stage order-four ROS method derived in \cite[Eqs. (35)-(37), pp. 153--154]{GonzalezPintoJCAM2017}.
	\end{itemize}
	Such ROS methods are A-stable and require solving some linear systems at each step whose coefficient matrices involve the exact Jacobian of the problem; otherwise, choosing ‘inexact Jacobian’, they would become W-methods, for which the order of convergence is generally lower, and furthermore there is no study on the related stability properties. Thanks to the results obtained in this paper, for TASE-RK methods we can instead consider ‘inexact Jacobian’ in the underlying linear systems to solve, and prove their unconditional stability even with this choice. For the implementation of TASE-RK, we have set as underlying explicit RK methods those reported in \cite[Appendix]{Bassenne2021}.
	
	\subsection{Burgers' equation}\label{subsecBurg}
	We consider the Burgers' equation \cite[Subsec. 3.2]{Calvo2021}
	\begin{equation}\label{Burg}
		u_t=\epsilon u_{xx}-\frac{1}{2} (u^2)_x, \quad \text{where } (x,t) \in [0,2\pi]\times[0,4],
	\end{equation}
    equipped with periodic Boundary Conditions (BCs) and Initial Condition (IC) $u_0=u(x,0)=(1-\cos(x))/2$. Fixing the grid $\{x_m=mh;\ m=0,\ldots,M;\ x_M=2\pi\}$, with $M=1024$, we discretize \eqref{Burg} in space using finite differences of order four, obtaining
    \begin{equation}\label{Burgsemidisc}
    	\vec{u}_t= \epsilon \mathcal{L}_1 \vec{u}-\frac{1}{2} \mathcal{L}_2 \vec{u}^2, \quad \vec{u}\in \mathbb{R}^M, \ \mathcal{L}_1, \ \mathcal{L}_2 \in \mathbb{R}^{M \times M}, \quad \text{where}
    	    \end{equation}
     $$\mathcal{L}_1=\frac{1}{h^2}\text{pentadiag}\bigg(-\frac{1}{12}, \frac{4}{3},-\frac{5}{2},\frac{4}{3},-\frac{1}{12}\bigg), \ \mathcal{L}_2=\frac{1}{h}\text{pentadiag}\bigg(\frac{1}{12},-\frac{2}{3},0,	\frac{2}{3},-\frac{1}{12}\bigg).$$
    Note that, since we are using periodic BCs, the elements in the corners of the $\mathcal{L}_1$ and $\mathcal{L}_2$ matrices are not zero. In this example, we apply TASE-RK methods \eqref{TASE_RK} setting
    \begin{equation}\label{BurgA}
    	A=\kappa (\epsilon \mathcal{L}_1 - 2I) \in \mathbb{R}^{M\times M}, \quad \kappa>0.
    \end{equation}
    That is, we perturb the diagonal elements of $\epsilon \mathcal{L}_1$ (in order to take $A$ symmetric negative definite), and introduce the parameter $\kappa$ according to Remark \ref{rem5}. To apply the theory developed in this manuscript, we set $B:=J_0-A$, where $J_0$ is the Jacobian of \eqref{Burgsemidisc} evaluated at $\vec{u}_0$ (see Remark \ref{rem4}). The matrices $A$ and $B$ do not commute.

\begin{figure}[h!]
    	\centering
    	\includegraphics[scale=0.435]{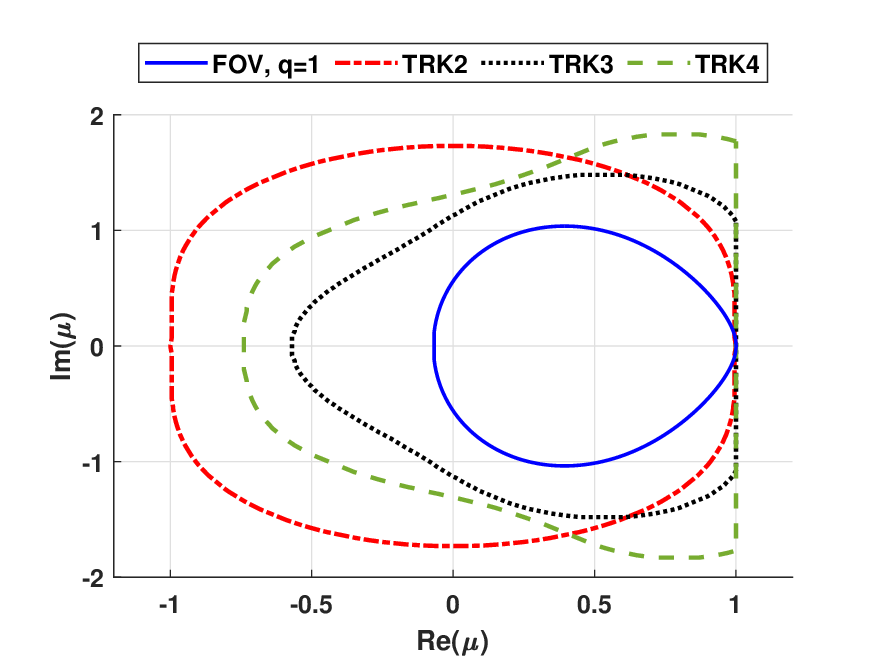}
    	\includegraphics[scale=0.435]{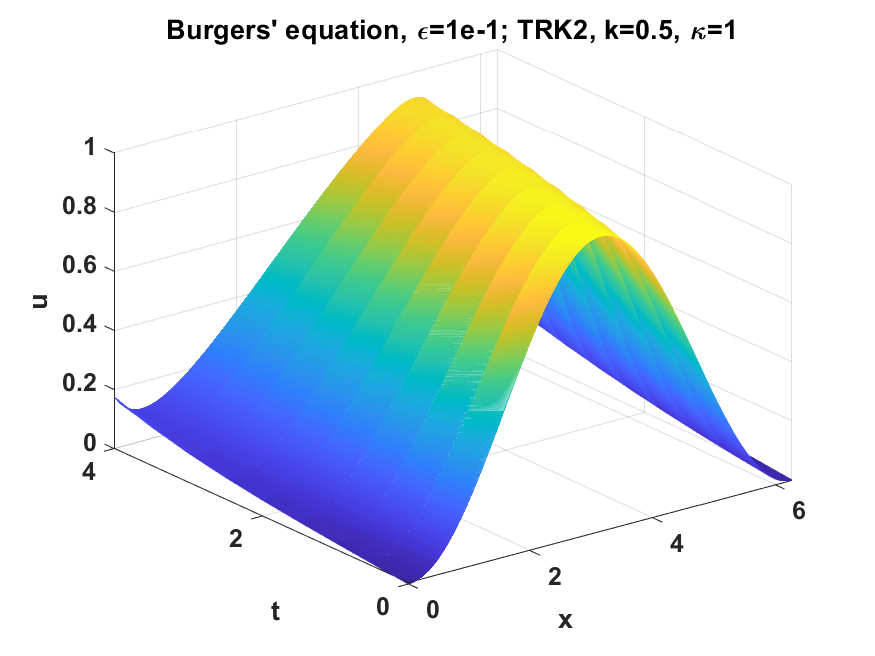}
    			\includegraphics[scale=0.435]{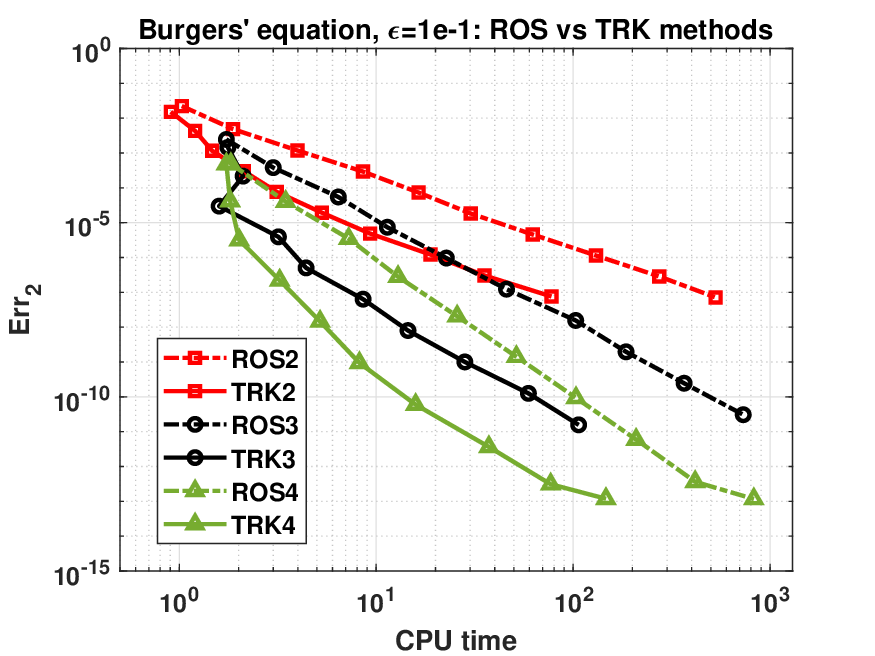}
    	\includegraphics[scale=0.435]{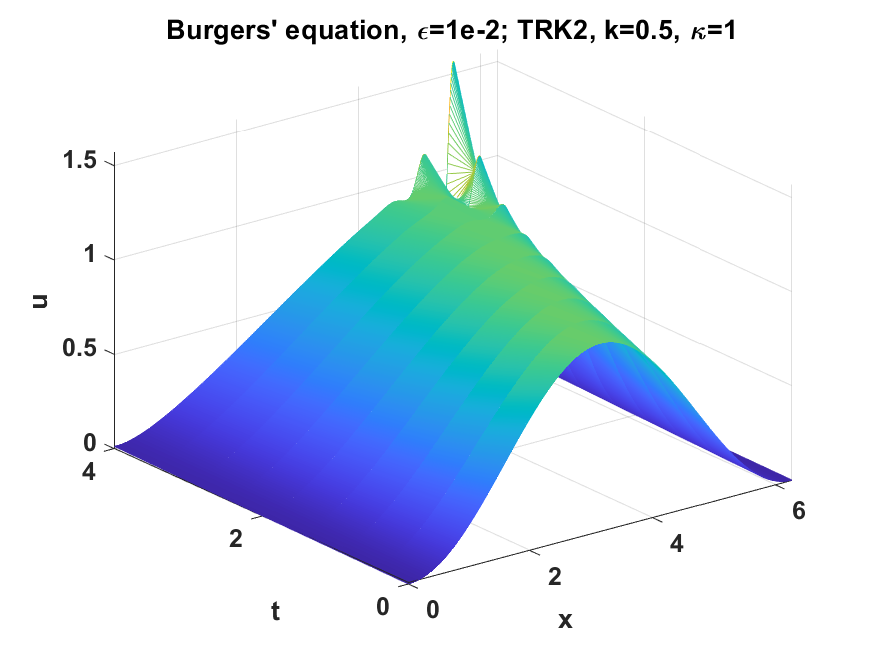}
    	\caption{$\mathcal{W}_{1}(-A,B)$, with $A$ as in \eqref{BurgA}, $\epsilon=10^{-1}, \ \kappa=1$, and $-\mathcal{ D}_{\infty,p}$ for $p=2,3,4$ (top left); solution of Burgers' equation by TRK2 method, with $k=0.5, \ \kappa=1$, and $\epsilon=10^{-1}$ (top right), $\epsilon=10^{-2}$ (bottom right); comparison in terms of CPU time and error between ROS and TRK methods on Burgers' equation, $\epsilon=10^{-1}$ (bottom left).}\label{Burg1}
    \end{figure}

    We first solve \eqref{Burgsemidisc} by considering $\epsilon= 10^{-1}$. In this case, for $\kappa=1$ all the TASE-RK methods are unconditionally stable, according to Theorem \ref{th5.4}. Indeed, as shown in Figure \ref{Burg1} (top, left), $W_1(-A,B)$, with $\kappa=1$, is included in $-\mathcal{ D}_{\infty,p}$ for $p=2,3,4$. As an example of this, note that the numerical solution given by the TRK2 method with step-size $k=0.5$ has a stable trend, see Figure \ref{Burg1} (top, right). Figure \ref{Burg1} (bottom, left) shows a comparison between TASE-RK and ROS methods, in terms of error and CPU time: the relative error (in norm two) is computed at the last grid point $t_e=4$ with respect to a reference solution given by the Matlab routine \texttt{ode15s}, setting \texttt{AbsTol=RelTol=eps}; the tests have been performed using the version R2023b of Matlab on a laptop with a RAM of 8GB and an operative system of 64 bit, where the processor is AMD Ryzen 7 3700U with Radeon Vega Mobile Gfx 2.30 GHz. Note that the TASE-RK are much more efficient than the ROS of the same order.

\begin{figure}[h!]
		\centering
		\includegraphics[scale=0.435]{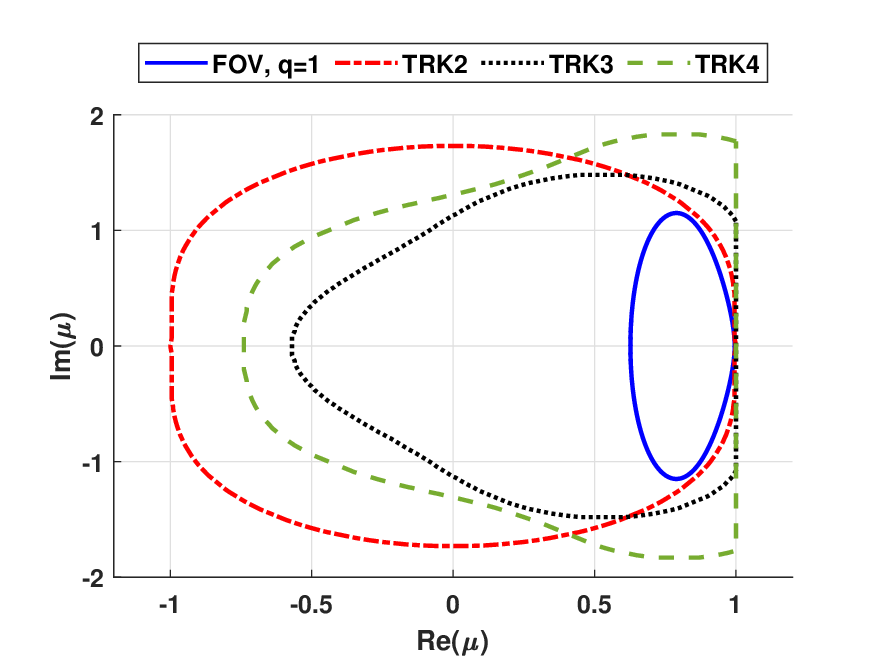}
		\includegraphics[scale=0.435]{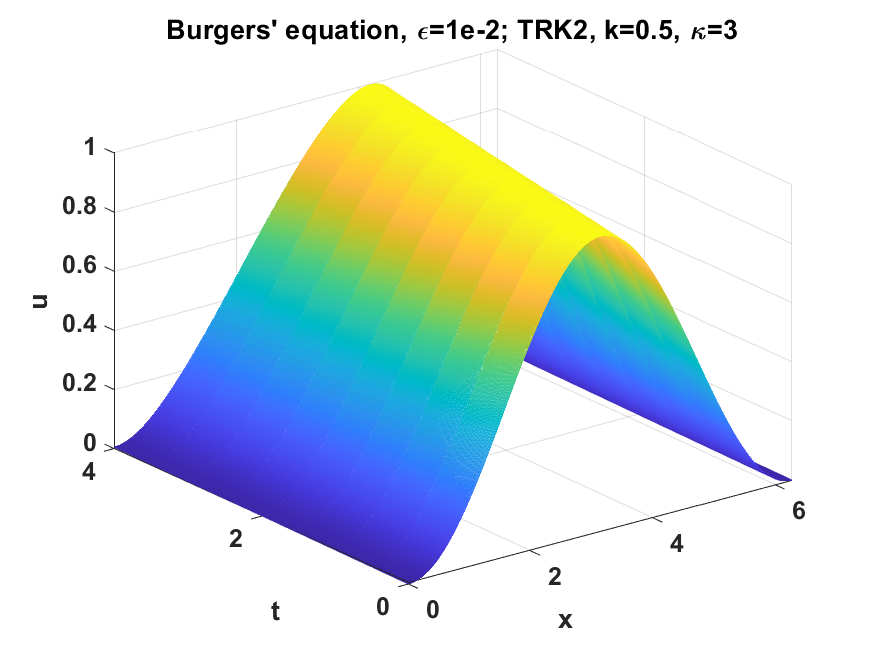}
		\caption{$\mathcal{W}_{1}(-A,B)$, with $A$ as in \eqref{BurgA}, $\epsilon=10^{-2}, \ \kappa=3$, and $-\mathcal{ D}_{\infty,p}$ for $p=2,3,4$ (left); solution of Burgers' equation by TRK2 method, with $\epsilon=10^{-2}, \ k=0.5, \ \kappa=3$ (right).}\label{Burg3}
	\end{figure}

	Finally, we also consider the Equation \eqref{Burgsemidisc} with $\epsilon=10^{-2}$. In this case, for $\kappa=1$ in \eqref{BurgA}, the TASE-RK methods are not unconditionally stable, since e.g. $-0.4875 \pm 3.4130i$ ($\in \upmu(A,B)$) $\notin \mathcal{ D }_{\infty,p}$ (see the necessary condition given by Theorem \ref{th5.4}). As an example to confirm this, see Figure \ref{Burg1} (bottom, right): the TRK2 displays instabilities for $k=0.5$. Setting instead $\kappa=3$ in \eqref{BurgA}, using Theorem \ref{th5.4}, all the TASE-RK methods are unconditionally stable, see Figure \ref{Burg3} (left). Indeed, with $\kappa=3$, now the TRK2 has a stable trend for $k=0.5$, see Figure \ref{Burg3} (right).

\subsection{FitzHugh-Nagumo model}
We consider the following FitzHugh-Nagumo (FHN) model \cite[Ch. 30]{bookMatlabFHN}, that describes the reaction of a biological neuron when excited by an external stimulus:
\begin{equation}\label{FHN}
	\begin{cases}
		\displaystyle v_t=D  v_{xx}+v-\frac{v^3}{3}-w,\\
		\displaystyle w_t=\frac{1}{\tau} (v-a-bw),
	\end{cases}\quad (x,t) \in [0,10]\times[0,t_e].
\end{equation}
To reproduce two repeating waves, we set $D=0.01$, $a=-0.7$, $b=0.8$, $\tau=12.5$, considering periodic BCs and, with $\phi(x):=1 + e^{-10\pi \sinh(x)}$, IC
\begin{equation*}
			\displaystyle v(x,0)= -\frac{3}{2} + \frac{3}{\phi(x-1.5)} - \frac{3}{\phi(x-2)} , \quad w(x,0)=\frac{-3}{4\phi(x-1.5)}.
\end{equation*}

As for Burgers' equation we perform a spatial discretization with $M=1024$ points, using finite differences of order four for $v_{xx}$; thus, the first equation in \eqref{FHN} turns into $ \vec{v}_t=D\mathcal{L}_1 \vec{v}+\vec{v}-\frac{\vec{v}^3}{3}-\vec{w}$, where $\mathcal{L}_1$ is as in \eqref{Burgsemidisc}. To apply TASE-RK methods, we set
\begin{equation}\label{FHNA}
A= \kappa \cdot \text{blkdiag} (D \mathcal{L}_1 - I, -(b/\tau) I) \in \mathbb{R}^{2M\times 2M}, \quad \kappa>0.
\end{equation}
In this way, $A$ is symmetric negative definite. To apply the theory of this paper, we define $B:=J_0-A$, where $J_0$ is the Jacobian of the FHN model at $\vec{u}_0$.

\begin{figure}[h!]
	\centering
	\includegraphics[scale=0.435]{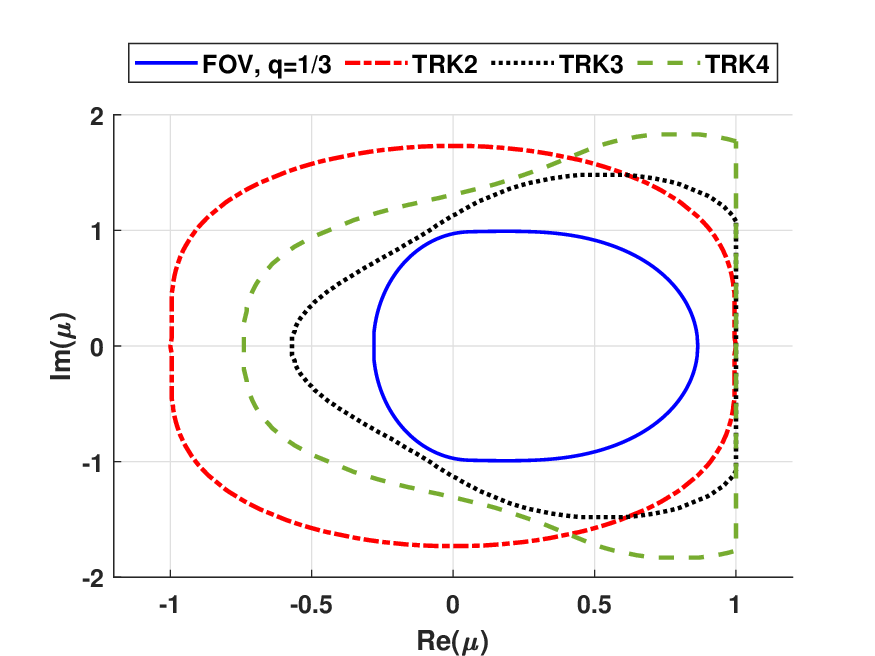}
	\includegraphics[scale=0.435]{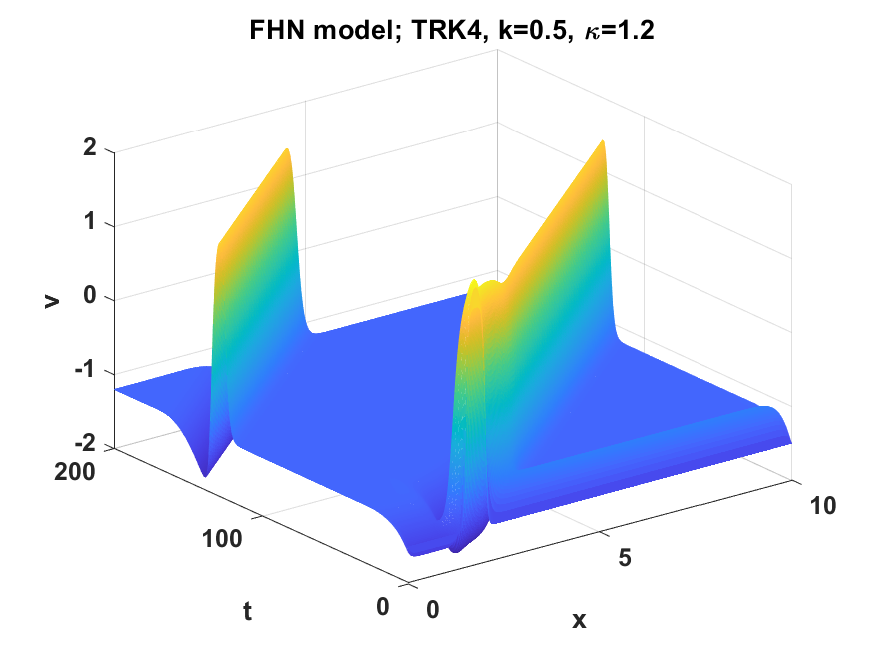}
	\includegraphics[scale=0.435]{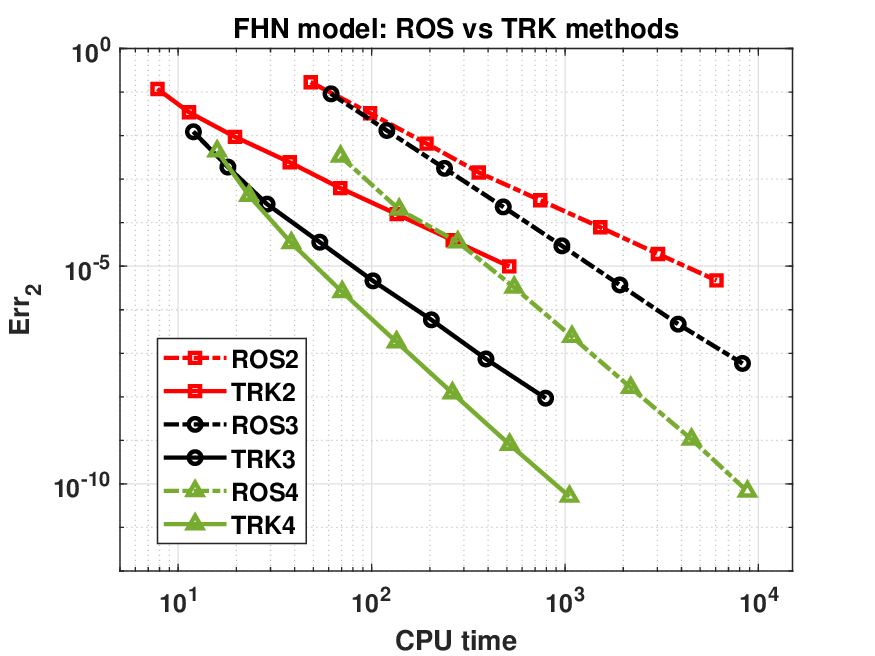}
	\caption{$\mathcal{W}_{1/3}(-A,B)$, with $A$ as in \eqref{FHNA}, $\kappa=1.2$, and $-\mathcal{ D}_{\infty,p}$ for $p=2,3,4$ (top, left); solution of FHN model by TRK4 method, with $k=0.5, \ \kappa=1.2, \ t \in [0,200]$ (top, right); comparison between ROS and TRK methods on FHN model, with $\kappa=1.2, \ t_e=50$ (bottom).}\label{FHNFig}
\end{figure}

Since $A$ and $B$ do not commute, we need to use the FOV to investigate the stability of TASE-RK methods. In particular, using Theorem \ref{th5.4}, for $\kappa=1.2$ in \eqref{FHNA}, all TASE-RK methods are unconditionally stable, as can be seen in Figure \ref{FHNFig} (top, left). In Figure \ref{FHNFig} (top, right) we e.g. report the numerical solution $v$ given by the TRK4 for $k=0.5, \ \kappa=1.2$, in the time interval $[0,200]$; as expected, the solution has a stable trend. Finally, we report in Figure \ref{FHNFig} (bottom) a comparison between TASE-RK and ROS methods (the tests and the plots were carried out as explained in Subsection \ref{subsecBurg}), which testifies the better efficiency of the former.

\section{Conclusions and future perspectives} \label{sec8}
In this manuscript, we have studied the stability properties of TASE-RK schemes with inexact Jacobian $J$, through an approach that is generally used for ImEx methods. In particular, we have considered a Jacobian splitting of the type $J=A+B$, analyzing both the case in which the matrices $A$ and $B$ are simultaneously diagonalizable and not. The numerical experiments, conducted on PDEs of interest in applications, have highlighted the correctness of the theoretical analysis, the good stability properties of the TASE-RK methods with inexact Jacobian, and their better efficiency compared to some ROS methods.

Several interesting works can be done starting from this paper. First of all, some issues raised in the manuscript could be further explored, e.g. the case in which $A$ and $B$ are simultaneously block-diagonalizable (Test \ref{test2}), and nonlinear stability (Section \ref{sec6}). Furthermore, one could try to find a setting of the free parameters $\omega_j, \ j=1,\ldots,p$, of the TASE operators \eqref{TASE_Montijano} that allows maximizing the area of the unconditional stability diagrams $\mathcal{ D }_{\infty,p}$, for $p=2,3,4$. Finally, it would be really interesting to use the approach proposed in the manuscript to analyze the stability of general W-methods (recall that TASE-RK are a subclass of W-methods).

\appendix
\section{Proofs of Proposition \ref{prop1}, Lemma \ref{lemmaC}, Proposition \ref{prop6.1}}
	
\subsection{Proposition \ref{prop1}}	
	\begin{proof}
		Remembering Definition \ref{def1}, let us take into account the test IVP \eqref{EqTest1mat} with $\vec{g}={\bf{0}}$. For the simultaneous diagonalizability of $A$ and $B$, it holds that
		\begin{equation*}
			A=Q \Lambda Q^{-1}, \quad \text{with} \quad \Lambda=\text{diag}(\lambda_i), \qquad B = Q \Gamma Q^{-1}, \quad \text{with} \quad \Gamma=\text{diag}(\gamma_i).
		\end{equation*}
		Setting $\vec{u}=Q\vec{\tilde{u}}$, the considered test problem turns into
		$
			\vec{\tilde{u}}_t=\Lambda\vec{\tilde{u}} +  \Gamma \vec{\tilde{u}}$.
		We can then write it as a system of decoupled scalar equations of the form
		$(\vec{\tilde{u}}_t)_i=\lambda_i(\vec{\tilde{u}})_i+\gamma_i(\vec{\tilde{u}})_i$, $i=1,\ldots,d$. Applying to them a TASE-RK method \eqref{TASE_RK} with $s=p\leq 4$ leads to
		\begin{equation}\label{recurrence}
		{(\vec{\tilde{u}}_{n+1})}_i= (\tilde{R}T_p(y_i,\mu_i))^{n+1} {(\vec{\tilde{u}}_{n})}_i,
		\end{equation}
		where $\tilde{R}T_p$ is as in \eqref{Rtilde}.
		Therefore, considering for each $y_i$ the related stability diagram $\mathcal{D}_{y_i,p}=\{\mu_i \in \mathbb{C}: |\tilde{R}T_p(y_i,\mu_i)|\leq1,  \text{ with Re}(\mu_i)\geq-1 \}$, we get the proof.
	\end{proof}
	
\subsection{Lemma \ref{lemmaC}}
	\begin{proof}
		There exist two matrices $P$ and $D$ such that $X =PDP^{-1}$: if $X$ is diagonalizable, $D$ is diagonal; otherwise, $D$ is the so-called Jordan normal form of $X$ \cite[Theorem 3.1.11]{HornJohnsonMATAN}. For greater generality, let us suppose we are in the latter case. Therefore, the matrix $D$ can be expressed as
		\begin{equation*}	
			D=\text{blkdiag}(D_{m_j}(c_j))_{j=1}^{\bar{d}},
			\ \sum_{j=1}^{\bar{d}}m_j=d, \ \ \ D_{m_j}(c_j)=\begin{pmatrix}
				c_j & 1 & & \\
				& \ddots & \ddots & \\
				&  & c_j & 1 \\
				&  & & c_j \\
			\end{pmatrix} \in \mathbb{C}^{m_j \times m_j}.
		\end{equation*}
		Note that $c_j$, $j=1,\ldots,\bar{d}$, constitute all the eigenvalues of $D$, and thus of $X$ (some of them could have multiplicity greater than one, and for this we have used the index $\bar{d}\leq d)$.	
		Since $D=D_{m_1}(c_1)\oplus \ldots \oplus D_{m_{\bar{d}}}(c_{\bar{d}})$ ($\oplus$ stands for the direct sum of matrices), we have that
		$
			D^q=D_{m_1}^q(c_1)\oplus \ldots \oplus D_{m_{\bar{d}}}^q(c_{\bar{d}}), \ \forall q\geq 1$,
		with \cite[Subsec. 3.2.5]{HornJohnsonMATAN}
		\begin{equation*}	
			D^q_{m_j}(c_j)=\begin{pmatrix}
				c_j^q & \binom{q}{1} c_j^{q-1} & \binom{q}{2} c_j^{q-2} & \dots & \binom{q}{m_j-1} c_j^{q-m_j+1} \\
				& c_j^q & \binom{q}{1} c_j^{q-1} & \dots & \binom{q}{m_j-2} c_j^{q-m_j+2} \\
				&  & \ddots & & \vdots \\
				&  & & c_j^q & \binom{q}{1} c_j^{q-1} \\
				&  & &  & c_j^q \\
			\end{pmatrix} \in \mathbb{C}^{m_j \times m_j}, \ j=1,\ldots,\bar{d}.
		\end{equation*}
		Thus, $ \sigma(D^q)=\{c_i^q, \ i=1,\ldots,d \}$ (we reuse the index $d$, so we count $m$ times an eigenvalue with multiplicity $m$).
		Since $I,D,D^2,\ldots,D^p$ are upper triangular matrices whose i-th diagonal elements correspond to $1,c_i,c_i^2,\ldots,c_i^p$, respectively, we have
		\begin{equation*}
			(I+D+D^2+\ldots+D^p)_{ij \ (j\leq i)}=\begin{cases} \begin{split}  & 1+c_i+c_i^2+\ldots + c_i^p, & \text{ if } j =i, & \\ & 0, & \text{ if } j<i, \end{split} \end{cases}	
			\text{ and therefore}
		\end{equation*}
		\begin{equation}\label{rel22}
		\sigma(\sum_{q=0}^pw_qD^q ) = \{\sum_{q=0}^pw_jc_i^q, \ i=1,\ldots,d \}.
		\end{equation}
		Remembering $X=PDP^{-1}$, and noting that $(PDP^{-1})^q=PD^qP^{-1}\ \forall q\geq 1$, we have
		\begin{equation}\label{rel3}
		\sigma(\sum_{q=0}^pw_qX^q )=\sigma(P(\sum_{q=0}^pw_qD^q)P^{-1} )=\sigma(\sum_{q=0}^pw_qD^q).
		\end{equation}
		From \eqref{rel22} and \eqref{rel3}, the theorem follows.
	\end{proof}
	
\subsection{Proposition \ref{prop6.1}}
	\begin{proof}
		Recall, from \eqref{FOVexpr}, that $\mathcal{W}_1(-A,B)=\{ \langle \vec{v}, B \vec{v} \rangle : \langle \vec{v}, -A \vec{v} \rangle = 1 \}$, and $\mathcal{ W }(-A)=\{ \langle \vec{w}, -A \vec{w} \rangle : ||\vec{w}||^2=1 \}$, with $\vec{v}, \ \vec{w} \in \mathbb{C}^d$. Condition $\langle \vec{v}, -A \vec{v} \rangle = 1$ can be equivalently expressed as
		$< \frac{\vec{v}}{||\vec{v}||}, -A \frac{\vec{v}}{||\vec{v}||}> = \frac{1}{||\vec{v}||^2}$, i.e. $\langle \vec{w}, -A \vec{w} \rangle = \frac{1}{||\vec{v}||^2}$, with $\vec{w}:=\frac{\vec{v}}{||\vec{v}||}$.
		Observing that $||\vec{w}||^2=1$, we therefore get
		\begin{equation}\label{W1}
		\mathcal{W}_1(-A,B)=\{ \langle \vec{v}, B \vec{v} \rangle : 1/\langle \vec{v}, \vec{v} \rangle \in \mathcal{ W }(-A) \}.
		\end{equation}
		Setting $\ell:=1/\langle \vec{v},\vec{v} \rangle$, and then $\vec{y}:=\ell ^{1/2} \vec{v}$, from \eqref{W1},
		\begin{equation}\label{W12}
		\mathcal{W}_1(-A,B)=\bigg\{ \biggl< \vec{y},  \frac{1}{\ell} B \vec{y} \biggr> : \langle \vec{y}, \vec{y} \rangle = 1, \ \ell \in \mathcal{ W }(-A) \bigg\}.
		\end{equation}
		From \eqref{W12}, using that $B=B(\xi)=\epsilon \ \text{diag}(1-2\xi_i)_{i=1}^{M-1}$, we can compute $\mathcal{W}_1(-A,B)$ and related bounds through the same steps as in \cite[Appendix A]{Seibold2019}, leading to
		\begin{equation*}
			\mathcal{W}_1(-A,B)=\frac{\epsilon}{\ell}   (1-2 \sum_{i=1}^{M-1} \xi_i|y_i|^2), \qquad \ell \in \mathcal{ W }(-A), \ 0 \leq |y_i|^2 \leq 1.
		\end{equation*}
		Let us consider two bounds (lower and upper) for $\xi_i$, i.e. $\upsilon_{\text{LB}}\leq \xi_i \leq \upsilon_{\text{UB}} \ \forall i=1,\dots,M-1$. Therefore, using that $\epsilon>0$, and observing that $\ell>0$,
		\begin{equation}\label{A6}
		\frac{\epsilon}{\ell}	(1-2 \upsilon_{\text{UB}}) \leq \mathcal{W}_1(-A,B) \leq    \frac{\epsilon}{\ell}	(1-2 \upsilon_{\text{LB}}), \qquad \ell \in \mathcal{ W }(-A).
		\end{equation}
		
		Finally, we need to characterize the elements of $\mathcal{ W }(-A)$. Since $-A$ is a symmetric positive definite matrix, from Remark \ref{remFOV}, item 5, it holds $\mathcal{ W }(-A)=\text{co}(\sigma(-A))$, and furthermore $\sigma(-A)$ takes values in $\mathbb{R}^+$. In particular, with $A$ as in \eqref{A},
		\begin{equation}\label{A7}
		\sigma(-A)=\bigg\{\frac{2D}{h^2} \bigg(1 + \cos\Big(\frac{i \pi}{M}\Big)\bigg), \ i=1,\ldots,M-1\bigg\}.
		\end{equation}
		Considering \eqref{A6} and \eqref{A7}, we get the proof.
	\end{proof}

\section*{Acknowledgments}
	The authors Conte, Pagano and Paternoster are members of the GNCS group. This work has been supported by: GNCS-INdAM project; Italian Ministry of University and Research (MUR), through PRIN PNRR 2022 project No. P20228C2PP and PRIN 2020 project No. 2020JLWP23; %“Integrated Mathematical Approaches to Socio–Epidemiological Dynamics” %(CUP: E15F21005420006)
	“Visiting Professors mobility Program” of the University of Salerno; %;
	Spanish Ministerio de Ciencia e Innovación (MCIN) with funding from the European Union NextGenerationEU (PRTRC17.I1); Consejer\'ia de Educación, Junta de Castilla y Le\'on, through QCAYLE project; Fundaci\'on Sol\'orzano through FS/5-2022 project.  	
	The authors thank Prof. D. Shirokoff (co-author of \cite{Rosales20172336,Seibold2019}), for conversations regarding: unconditional stability; stability of RK methods to PDEs; use of Lyapunov theory in the context of Remark \ref{rem4}.
	
	\bibliographystyle{siamplain}
	\bibliography{references}

\begin{thebibliography}{10}

\bibitem{Aceto2023}
{\sc L.~Aceto, D.~Conte, and G.~Pagano}, {\em On a generalization of
  time-accurate and highly-stable explicit operators for stiff problems}, Appl.
  Numer. Math., ISSN 0168-9274 (2023),
  \url{https://doi.org/10.1016/j.apnum.2023.04.001}.

\bibitem{Bassenne2021}
{\sc M.~Bassenne, L.~Fu, and A.~Mani}, {\em Time-accurate and highly-stable
  explicit operators for stiff differential equations}, J. Comput. Phys., 424
  (2021), \url{https://doi.org/10.1016/j.jcp.2020.109847}.

\bibitem{Boscarino2015}
{\sc S.~Boscarino, R.~B\"{u}rger, P.~Mulet, G.~Russo, and L.~M. Villada}, {\em
  Linearly implicit {IMEX} {R}unge--{K}utta methods for a class of degenerate
  convection-diffusion problems}, SIAM J. Sci. Comput., 37 (2015),
  pp.~B305--B331, \url{https://doi.org/10.1137/140967544}.

\bibitem{Butcherbook}
{\sc J.~C. Butcher}, {\em Numerical Methods for Ordinary Differential
  Equations}, John Wiley \& Sons, Ltd, 2008,
  \url{https://doi.org/10.1002/9780470753767}.

\bibitem{Calvo2023}
{\sc M.~Calvo, L.~Fu, J.~I. Montijano, and L.~Rández}, {\em Singly {TASE}
  operators for the numerical solution of stiff differential equations by
  explicit {R}unge–{K}utta schemes}, J. Sci. Comput., 96 (2023),
  \url{https://doi.org/10.1007/s10915-023-02232-3}.

\bibitem{Calvo2021}
{\sc M.~Calvo, J.~I. Montijano, and L.~Rández}, {\em A note on the stability
  of time–accurate and highly–stable explicit operators for stiff
  differential equations}, J. Comput. Phys., 436 (2021),
  \url{https://doi.org/10.1016/j.jcp.2021.110316}.

\bibitem{Conte2023}
{\sc D.~Conte, G.~Pagano, and B.~Paternoster}, {\em Nonstandard finite
  differences numerical methods for a vegetation reaction–diffusion model},
  J. Comput. Appl. Math., 419 (2023),
  \url{https://doi.org/10.1016/j.cam.2022.114790}.

\bibitem{Conte20231}
{\sc D.~Conte, G.~Pagano, and B.~Paternoster}, {\em Time-accurate and
  highly-stable explicit peer methods for stiff differential problems}, Comm.
  Nonlinear Sci. Numer. Simul., 119 (2023),
  \url{https://doi.org/10.1016/j.cnsns.2023.107136}.

\bibitem{chebfun}
{\sc T.~A. Driscoll, N.~Hale, and L.~N. Trefethen}, {\em Chebfun Guide},
  Pafnuty Publications, Oxford, 2014,
  \url{https://doi.org/www.chebfun.org/docs/guide/}.

\bibitem{Frasca-Caccia2023}
{\sc G.~Frasca-Caccia, C.~Valentino, F.~Colace, and D.~Conte}, {\em An overview
  of differential models for corrosion of cultural heritage artefacts}, Math.
  Model. Nat. Phenom., 18 (2023), \url{https://doi.org/10.1051/mmnp/2023031}.

\bibitem{GonzalezPintoJCAM2017}
{\sc S.~Gonz\'{a}lez-Pinto, D.~Hern\'{a}ndez-Abreu, and
  S.~P\'{e}rez-Rodr\'{\i}guez}, {\em {$W$}-methods to stabilize standard
  explicit {R}unge-{K}utta methods in the time integration of
  advection-diffusion-reaction {PDE}s}, J. Comput. Appl. Math., 316 (2017),
  pp.~143--160, \url{https://doi.org/10.1016/j.cam.2016.08.026}.

\bibitem{GonzalezPintoAMC2016}
{\sc S.~Gonz\'{a}lez-Pinto, D.~Hern\'{a}ndez-Abreu,
  S.~P\'{e}rez-Rodr\'{\i}guez, and R.~Weiner}, {\em A family of three-stage
  third order {AMF}-{W}-methods for the time integration of advection diffusion
  reaction {PDE}s}, Appl. Math. Comput., 274 (2016), pp.~565--584,
  \url{https://doi.org/10.1016/j.amc.2015.10.013}.

\bibitem{Gonzalez-Pinto2023129}
{\sc S.~González-Pinto, D.~Hernández-Abreu, G.~Pagano, and
  S.~Pérez-Rodríguez}, {\em Generalized {TASE}-{RK} methods for stiff
  problems}, Appl. Numer. Math., 188 (2023), pp.~129--145,
  \url{https://doi.org/10.1016/j.apnum.2023.03.007}.

\bibitem{HairerWanner}
{\sc E.~Hairer and G.~Wanner}, {\em Solving ordinary differential equations
  {II}}, vol.~14 of Springer Series in Computational Mathematics,
  Springer-Verlag, Berlin, second~ed., 1996,
  \url{https://doi.org/10.1007/978-3-642-05221-7}.

\bibitem{HornJohnsonFOV}
{\sc R.~A. Horn and C.~R. Johnson}, {\em Topics in Matrix Analysis}, Cambridge
  University Press, 1991, \url{https://doi.org/10.1017/CBO9780511840371}.

\bibitem{HornJohnsonMATAN}
{\sc R.~A. Horn and C.~R. Johnson}, {\em Matrix Analysis}, Cambridge University
  Press, second~ed., 2012, \url{https://doi.org/10.1017/CBO9780511810817}.

\bibitem{Hundsdorferbook}
{\sc W.~Hundsdorfer and J.~Verwer}, {\em Numerical Solution of Time-Dependent
  Advection-Diffusion-Reaction Equations}, vol.~33 of Springer Series in
  Computational Mathematics, Springer, 2003,
  \url{https://doi.org/10.1007/978-3-662-09017-6}.

\bibitem{JohnsonSINUM1978}
{\sc C.~R. Johnson}, {\em Numerical determination of the field of values of a
  general complex matrix}, SIAM J. Numer. Anal., 15 (1978), pp.~595--602,
  \url{https://doi.org/10.1137/0715039}.

\bibitem{LangNumerMath2013}
{\sc J.~Lang and J.~G. Verwer}, {\em W-methods in optimal control}, Numer.
  Math., 124 (2013), pp.~337--360,
  \url{https://doi.org/10.1007/s00211-013-0516-x}.

\bibitem{Murraybook}
{\sc J.~D. Murray}, {\em Mathematical Biology I, An Introduction},
  Interdisciplinary Applied Mathematics, Springer, New York, third~ed., 2002,
  \url{https://doi.org/doi.org/10.1007/b98868}.

\bibitem{Rosales20172336}
{\sc R.~R. Rosales, B.~Seibold, D.~Shirokoff, and D.~Zhou}, {\em Unconditional
  stability for multistep {IMEX} schemes: Theory}, SIAM J. Numer. Anal., 55
  (2017), pp.~2336--2360, \url{https://doi.org/10.1137/16M1094324}.

\bibitem{Seibold2019}
{\sc B.~Seibold, D.~Shirokoff, and D.~Zhou}, {\em Unconditional stability for
  multistep {IMEX} schemes: Practice}, J. Comput. Phys., 376 (2019),
  pp.~295--321, \url{https://doi.org/10.1016/j.jcp.2018.09.044}.

\bibitem{Steihaug}
{\sc T.~Steihaug and A.~Wolfbrandt}, {\em An attempt to avoid exact {J}acobian
  and nonlinear equations in the numerical solution of stiff differential
  equations}, Math. Comput., 33 (1979), pp.~521--534,
  \url{https://doi.org/10.2307/2006293}.

\bibitem{Trefethen}
{\sc L.~N. Trefethen and D.~Bau}, {\em Numerical Linear Algebra}, SIAM, 1997,
  \url{https://doi.org/10.1137/1.9780898719574.fm}.

\bibitem{bookMatlabFHN}
{\sc P.~Wallisch}, {\em MATLAB for Neuroscientists}, Academic Press, San Diego,
  second~ed., 2014, \url{https://doi.org/10.1016/B978-0-12-383836-0.00030-8}.

\end{thebibliography}
\end{document}